\PassOptionsToPackage{
    colorlinks=true,
    linkcolor=red,
    citecolor=blue,
    urlcolor=purple,
    pdfborder={0 0 0}
}{hyperref}

\documentclass[preprint,12pt,authoryear,nopreprintline]{elsarticle}
\usepackage[
    margin=1.15in
]{geometry}
\usepackage{amssymb}
\usepackage{amsmath}
\usepackage{booktabs}
\usepackage[normalem]{ulem}
\usepackage{theorem}
\usepackage{bbm}

\numberwithin{equation}{section} 
\newtheorem{teo}{Theorem} 
\newtheorem{cor}{Corollary} 
\newtheorem{example}{Example} 

\newtheorem{remark}{Remark}

\usepackage{xcolor}
\usepackage{hyperref}

\journal{Stochastic Processes and their Applications}

\begin{document}

\hypersetup{
    colorlinks=true,
    linkcolor=red,
    citecolor=blue,
    urlcolor=purple
}

\begin{frontmatter}

\title{Concentration of empirical entropy and transfer entropy for non-regular chains with unbounded memory}

\tnotetext[dedication]{\emph{In honor of Antonio Galves, whose mathematical insight lives on through generations of mathematicians.}}

\author[UFRB]{Jos\'{e} O. Santana}
\author[UFSCar]{Ricardo F. Ferreira}

\affiliation[UFRB]{organization={Department of Mathematics, Federal University of Reconcavo da Bahia}, city={Amargosa}, postcode={45300-000}, state={BA}, country={Brazil}
 }
 
\affiliation[UFSCar]{organization={Department of Statistics, Federal University of São Carlos}, 
      city={São Carlos}, postcode={13565-905}, state={SP}, country={Brazil}
 }

\begin{abstract}
We study plug-in estimators of the entropy rate and the transfer entropy rate for stationary finite-alphabet chains with unbounded memory. Under a coalescent coupling-from-the-past assumption, we establish non-asymptotic concentration inequalities around the expectations of the estimators. Under additional uniform non-nullness and exponential $\beta$-mixing assumptions, we control the corresponding biases and obtain concentration around the entropy and transfer entropy rates, together with almost-sure consistency. The results apply to chains governed by $g$-functions that need not be globally continuous.
\end{abstract}


\end{frontmatter}

\section{Introduction}
\label{sec_introduction}

Transfer entropy is an information-theoretic measure of directed dependence between two stochastic processes. More precisely, it quantifies how much knowledge of the past of one process reduces the uncertainty about the next state of another process, beyond what is already explained by the latter's own past. Introduced by \citet{schreiber2000measuring} and, independently, by \citet{paluvs2001synchronization}, transfer entropy has since found applications in several areas, including neuroscience \citep{vicente2011transfer, izzi2024identifying} and finance \citep{marschinski2002analysing, dimpfl2012using}.

In this paper, random information sources are modeled as discrete-time, time homogeneous stochastic processes taking values in finite alphabets. Our main object of interest is the \emph{transfer entropy rate}, which describes the asymptotic directed predictive dependence between two stochastic processes. Unlike symmetric information-theoretic quantities, transfer entropy distinguishes the direction of dependence and is therefore often interpreted as an operational notion of causality in a predictive sense. In this respect, it provides a complementary viewpoint to the classical Granger framework \citep{granger1969investigating}, whose standard formulation is based on linear predictive models, while transfer entropy is
not restricted to linear dependence.

A natural statistical question is whether the entropy rate and the transfer entropy rate can be consistently estimated and, beyond consistency, how the corresponding estimators fluctuate in finite samples. As we shall see, these two quantities are closely related, and so are their plug-in estimators. The main goal of this paper is to characterize the non-asymptotic stochastic behavior of these estimators. More precisely, we establish concentration inequalities around their expectations under the coupling-from-the-past assumption, and, under the additional uniform non-nullness and exponential $\beta$-mixing assumptions, concentration inequalities around the corresponding entropy and transfer entropy rates (Theorems~\ref{thm:concentration_entropy_mean}, \ref{thm:concentration_entropy_rate}, \ref{thm:concentration_transfer_entropy_mean} and \ref{thm:concentration_transfer_entropy_rate}). Our results apply to chains with unbounded memory governed by $g$-functions that are not necessarily continuous.

The probabilistic framework underlying these results is that of $g$-functions. A $g$-function specifies a one-step transition distribution as a function of the entire past. This framework includes finite-order Markov chains but, more importantly for the present work, also allows for stochastic chains with unbounded memory, classically studied under the name \emph{chains with complete connections}. Systematic investigations of such processes date back to
\citet{onicescu1935chaines} and \citet{doeblin1937chaines}. The terminology \emph{$g$-function} was introduced by \citet{keane1972strongly} and has since become standard, particularly in the ergodic-theory literature \citep{ledrappier1974principe, walters1975ruelle, hulse1991uniqueness, johansson2008square}.

The information-theoretic rates considered above are well defined under the stationarity of the underlying process. A stationary probability measure compatible with a $g$-function is commonly referred to as a \emph{$g$-measure}. Thus, in order to study the estimation of the entropy rate and the transfer entropy rate, one must work with a stationary probability measure compatible with the underlying $g$-function. Uniqueness is not required for the definition of these rates for a fixed stationary chain, but it ensures that the $g$-function determines a unique stationary law. Classical sufficient conditions for the existence and uniqueness of $g$-measures are typically formulated in terms of the continuity properties of the $g$-function \citep{harris1955chains, keane1972strongly}. More recently, several works \citep{desantis2012backward, gallo2013perfect, gallo2013non,oliveira2015stochastic, ferreira2020non, gallesco2025uniqueness} have investigated stationary compatible measures in settings where the underlying $g$-functions may be discontinuous. The framework considered in this paper belongs to this latter setting, i.e., we do not assume global continuity of the $g$-function.

The estimation of entropy rates from empirical block frequencies has a long history. For chains with infinite memory, \citet{gabrielli2003fluctuations} studied the fluctuations of empirical entropies based on blocks of logarithmically increasing length, establishing a central limit theorem for the conditional empirical entropy under exponential loss of memory. Subsequently, \citet{chazottes2005large} established large-deviation principles for conditional and non-conditional empirical entropies of $g$-measures. For one-dimensional Gibbs measures, \citet{chazottes2011concentration} derived concentration inequalities for empirical-frequency entropy estimators, including concentration of the conditional empirical entropy around the entropy rate. More recently, \citet{chazottes2023gaussian} obtained Gaussian concentration bounds for plug-in entropy estimators of stochastic chains with unbounded memory under summability conditions on the variations or oscillations of the transition kernel. The results obtained here concern a different regime: the underlying $g$-function need not be globally continuous, and concentration is instead obtained from a coalescent coupling-from-the-past construction, with uniform non-nullness and exponential $\beta$-mixing used to control the bias when concentration around the entropy rate is considered (see Theorems \ref{thm:concentration_entropy_mean} and \ref{thm:concentration_entropy_rate}). Related exponential inequalities in the statistical analysis of variable-length Markov chains were obtained by \citet{galves2008exponential}, who derived finite-sample bounds for empirical context-tree estimation and established almost-sure consistency of the Context algorithm.

Several approaches have been proposed for estimating transfer entropy from observed time series. These include symbolic methods \citep{staniek2008symbolic}, nonparametric procedures based on discretization, kernel-density estimation, and adaptive partitioning \citep{verdes2005assessing, lee2012transfer}, bias-compensated estimators designed for finite and high-dimensional samples \citep{faes2013compensated}, and multiscale approaches
\citep{lungarella2007information}. For the transfer entropy rate specifically, \citet{restrepo2020transfer} proposed an estimator based on Lempel--Ziv complexity. These works provide a variety of computational approaches to transfer-entropy estimation, but comparatively less is known about the rigorous finite-sample fluctuations of the corresponding estimators. From a theoretical perspective, the estimation of directed information provides a closely related precedent. \citet{kontoyiannis2016estimating} studied the plug-in estimator of the directed information rate for finite-order Markov processes, establishing its asymptotic distribution and its connection with likelihood-ratio testing for causal influence. Building on this statistical framework, \citet{izzi2024identifying} considered the plug-in estimator of the transfer entropy rate for variable-length Markov chains with bounded memory, derived its asymptotic distribution under the null hypothesis of vanishing transfer entropy, and applied the resulting methodology to the identification of effective connectivity in stochastic neuronal networks. Moving from asymptotic distribution theory to non-asymptotic concentration, \citet{ferreira2021concentration} established concentration bounds for the plug-in estimator of directed information for chains with unbounded memory under continuity and exponential loss-of-memory conditions. In the present work, we study instead the plug-in estimator of the transfer entropy rate for stationary finite-alphabet chains with unbounded memory without requiring the underlying $g$-function to be globally continuous. We derive non-asymptotic concentration inequalities around its expectation under a coalescent coupling-from-the-past assumption, and, under additional uniform non-nullness and exponential $\beta$-mixing assumptions, obtain concentration around the transfer entropy rate itself (see Theorems \ref{thm:concentration_transfer_entropy_mean} and \ref{thm:concentration_transfer_entropy_rate}).

The remainder of the paper is organized as follows. In Section~\ref{sec:preliminary_notions}, we introduce the notation and preliminary notions used throughout the paper, with particular emphasis on chains with unbounded memory, entropy rates, transfer entropy rates, and their plug-in estimators. Section~\ref{sec:main_results} presents our main concentration results. In Section~\ref{sec:simulation}, we illustrate these results through numerical simulations for the non-Markovian models introduced in Section~\ref{sec:preliminary_notions}. The proofs of the main results are given in Section~\ref{sec:proofs}. Finally, Section~\ref{sec:final_remarks} contains some concluding remarks.

\section{Notations, definitions and preliminary notions}
\label{sec:preliminary_notions}

We begin by establishing the notation and basic concepts that will be used throughout the article.

\subsection{Notations}
\label{sesc:notations}

In this work, we consider a finite set $\mathcal{Z}$ endowed with the discrete topology, which we refer to as the \emph{alphabet}. We denote by $\mathbb{Z}:=\{\ldots,-1,0,1,\ldots\}$ the set of integers and take $\mathbb{N}:=\{0,1,2,\ldots\}$ as the set of natural numbers. We denote by $z:=(\ldots,z_{-1},z_0,z_1,\ldots)$ a generic element of $\mathcal{Z}^{\mathbb{Z}}$, endowed with the product topology; thus, $z\in\mathcal{Z}^{\mathbb{Z}}$ is a bi-infinite sequence of symbols from $\mathcal{Z}$. Given integers $j\leq k$, we write $[[j,k]]:=\mathbb{Z}\cap[j,k]$ and denote by $z_j^k:=(z_j,z_{j+1},\ldots,z_k)\in\mathcal{Z}^{[[j,k]]}$ the corresponding finite block, whose length is $|z_j^k|:=k-j+1$. By convention, if $j>k$, then $z_j^k$ denotes the empty block and $|z_j^k|=0$. Given $z_j^k\in\mathcal{Z}^{[[j,k]]}$, we denote by $[z_j^k]$ the associated cylinder set, namely 

\[ 
[z_j^k] := \left\{ w\in\mathcal{Z}^{\mathbb{Z}}: w_i=z_i \text{ for all } i\in[[j,k]] \right\}. \]

When $j=1$, we write $\mathcal{Z}^{[[1,k]]}$ simply as $\mathcal{Z}^k$. In particular, a block $z_{-k}^{-1}$ of length $k$ is naturally regarded as
an element of $\mathcal{Z}^k$. We denote by
\[
\mathcal{Z}^{*}:=\bigcup_{m\in\mathbb{N}}\mathcal{Z}^m,
\qquad
\mathcal{Z}^0:=\{\varnothing\},
\]
the set of all finite blocks over $\mathcal{Z}$, including the empty block. For every $k\in\mathbb{Z}$, we denote by $\mathcal{Z}^{\leq k}$ the set of left-infinite pasts $z_{-\infty}^k:=(\ldots,z_{k-1},z_k)$, and by $\mathcal{Z}^{\geq k}$ the set of right-infinite futures $z_k^{+\infty}:=(z_k,z_{k+1},\ldots).$

Given two finite blocks $z_j^k$ and $w_s^t$, with $j\leq k$ and $s\leq t$, their concatenation is denoted by $z_j^k w_s^t:=(z_j,\ldots,z_k,w_s,\ldots,w_t),$ and concatenation extends to semi-infinite sequences in the obvious way.

We denote random variables by uppercase letters, stochastic chains by bold uppercase letters, and their realizations by lowercase letters. In particular, $\boldsymbol{Z}:=\{Z_t:t\in\mathbb{Z}\}$ denotes a stochastic chain, and subscripts indicate time/position in a sequence.

We denote by $\mathcal{F}:=\mathcal{B}\bigl(\mathcal{Z}^{\mathbb{Z}}\bigr)$ the Borel $\sigma$-algebra of $\mathcal{Z}^{\mathbb{Z}}$, which, under the product topology, is equivalently the $\sigma$-algebra generated by the cylinder sets. Let $\pi_j:\mathcal{Z}^{\mathbb{Z}}\to\mathcal{Z}$, $\pi_j(z):=z_j$, denote the $j$-th coordinate projection. For each $k\in\mathbb{Z}$, we set $\mathcal{F}^{\leq k}:=\sigma(\pi_j:j\leq k)$ and
$\mathcal{F}^{\geq k}:=\sigma(\pi_j:j\geq k).$ Thus, $\mathcal{F}^{\leq k}$ contains the events determined by the coordinates up to time $k$, whereas $\mathcal{F}^{\geq k}$ contains those
determined by the coordinates from time $k$ onward.

We set $\mathbb{N}^{*}:=\mathbb{N}\setminus\{0\}$, $\overline{\mathbb{N}}:=\mathbb{N}\cup\{+\infty\}$ and $\overline{\mathbb{N}}^{*} = \mathbb{N}^{*} \cup\{+\infty\}$. Finally, $\mathbf{1}\{\cdot\}$ stands for the indicator function.

\subsection{Chains with unbounded memory}
\label{sec:chains_undounded_memory}

Chains with unbounded memory represent a natural generalization of finite-order Markov chains. In this setting, the probability of the next symbol may depend on the entire past of the process, in contrast to Markov chains in which the dependence is restricted to a fixed number of previous states.

The dynamics of such a process are described by a \emph{$g$-function}, that is, a family of transition probabilities specifying the law of the next symbol as a function of the infinite past. Formally, a \emph{$g$-function} is a measurable function $g:\mathcal{Z}^{\leq -1}\times\mathcal{Z}\longrightarrow[0,1]$
such that
\begin{align*}
    \forall \; z_{-\infty}^{-1} \in \mathcal{Z}^{\leq -1}, \quad 
    \sum_{z_0 \in \mathcal{Z}} g(z_0 \mid z_{-\infty}^{-1}) = 1.
\end{align*}
In other words, a $g$-function specifies a one-step transition distribution given the entire past, thus playing the role of a transition rule for chains with unbounded memory.

Since we are interested in time-homogeneous chains, we shall use the same $g$-function at every time index. With a slight abuse of notation, for each $i\in\mathbb{Z}$ and any $z\in\mathcal{Z}^{\mathbb{Z}}$, we write
\[
g(z_i \mid z_{-\infty}^{i-1})
:=
g\left(z_i \mid (T^{i}z)_{-\infty}^{-1}\right),
\]
where $T:\mathcal{Z}^{\mathbb{Z}}\rightarrow\mathcal{Z}^{\mathbb{Z}}$ denotes the \emph{left-shift operator}, that is,
\[
(Tz)_t=z_{t+1},
\qquad t\in\mathbb{Z}.
\]
More generally, for every $i,t\in\mathbb{Z}$,
\[
(T^i z)_t=z_{t+i}.
\]
In this way, the same transition rule is used at every time index, which reflects the homogeneity of the dynamics.

A stochastic chain $\boldsymbol{Z}:=\{Z_t:t\in\mathbb{Z}\}$ defined on a probability space $(\Omega,\Sigma,\mathbb{P})$, whose variables take values in the alphabet $\mathcal{Z}$, is said to be \emph{compatible} with a given $g$-function $g$ if, for every $i\in\mathbb{Z}$ and every $a\in\mathcal{Z}$,
\[
\mathbb{P}\left(
Z_i=a
\,\middle|\,
\sigma(Z_j:j\leq i-1)
\right)
=
g\left(a\mid Z_{-\infty}^{i-1}\right)
\qquad \mathbb{P}\text{-a.s.}
\]
Thus, the same conditional rule, prescribed by $g$, applies at every time index, which expresses the time homogeneity of the transition mechanism. Hence, chains compatible with a given $g$-function are governed by the transition probabilities prescribed by $g$. When these transition probabilities genuinely depend on arbitrarily remote coordinates of the past, the corresponding chain is said to have \emph{unbounded memory}. More precisely, the transition mechanism specified by $g$ is said to be of \emph{finite order} if there exists some $k\in\mathbb{N}^{*}$ such that, for every $a\in\mathcal{Z}$ and every $z_{-\infty}^{-1},w_{-\infty}^{-1}\in\mathcal{Z}^{\leq -1}$, \[ z_{-k}^{-1}=w_{-k}^{-1} \quad\Longrightarrow\quad g(a\mid z_{-\infty}^{-1}) = g(a\mid w_{-\infty}^{-1}). \] In this case, the transition probabilities depend only on the last $k$ symbols of the past. If no such finite $k$ exists, we say that the transition mechanism has \emph{unbounded memory}. In other words, under unbounded memory, the portion of the past needed to determine the transition probabilities is not uniformly bounded over all pasts. Depending on the past, it may be finite but arbitrarily long, or it may be genuinely infinite.

Let $\mu$ denote the law of $\boldsymbol{Z}$ on $(\mathcal{Z}^{\mathbb{Z}},\mathcal{F})$. Then, for every
$a\in\mathcal{Z}$,
\[
\mathbb{E}_{\mu}\left[
\mathbf{1}_{[a]}
\,\middle|\,
\mathcal{F}^{\leq -1}
\right](z)
=
g\left(a\mid z_{-\infty}^{-1}\right)
\qquad \mu\text{-a.s.}
\]
Although the law $\mu$ of a given chain is uniquely determined by the chain,
a given $g$-function may admit more than one compatible probability measure.

A probability measure $\mu$ on $(\mathcal{Z}^{\mathbb{Z}},\mathcal{F})$ is said to be \emph{shift-invariant} if \[ \mu(B)=\mu(T^{-1}B), \qquad \forall\, B\in\mathcal{F}. \] A shift-invariant probability measure $\mu$ that is compatible with the $g$-function $g$ is called a \emph{$g$-measure}. Equivalently, under a $g$-measure, the canonical stochastic chain is stationary and its transition probabilities are prescribed by $g$. The existence of a $g$-measure, however, is not automatic for a general $g$-function.

Classical examples help illustrate how $g$-functions extend familiar models. In the i.i.d. case, the transition rule is independent of the past and is given by
\[
g(z_0 \mid z_{-\infty}^{-1})
=
\mathbb{P}(Z_0=z_0),
\qquad z_0\in\mathcal{Z}.
\]
More generally, for a homogeneous Markov chain of order $k\in\mathbb{N}^{*}$, the next symbol depends only on the last $k$ symbols of the past, and hence
\[
g(z \mid z_{-\infty}^{-1})
=
Q(z_{-k}^{-1},z),
\qquad z\in\mathcal{Z},
\]
where $Q$ is the corresponding transition kernel. Thus, i.i.d. sequences and finite-order Markov chains arise as particular cases of chains described by $g$-functions.

These examples also illustrate that, in classical settings, the associated $g$-measure is well understood. In the i.i.d. case, it exists and is unique. For a homogeneous Markov chain of order $k$ on a finite alphabet, a stationary measure exists, and it is unique under the usual
irreducibility assumption when the chain is viewed on the state space $\mathcal{Z}^k$. For general $g$-functions with unbounded memory, however, existence and uniqueness can no longer be
deduced directly from the finite-state Markov-chain framework and may require additional assumptions. This point is illustrated in the following example.

\begin{example}[A non-Markovian $g$-function]
\label{ex:non-Markovian-g-function}

Let the alphabet be $\mathcal{Z}=\{0,1\}$ and let $(q_n)_{n\in\overline{\mathbb{N}}^{*}}$
be a sequence of real numbers in $[0,1]$. We define the $g$-function $g:\{0,1\}^{\leq-1}\times\{0,1\}\rightarrow[0,1]$ by
\[
g(1\mid z_{-\infty}^{-1}) := q_{\ell(z_{-\infty}^{-1})},
\qquad
g(0\mid z_{-\infty}^{-1}) := 1-q_{\ell(z_{-\infty}^{-1})},
\]
where, for any
$z_{-\infty}^{-1}\in\mathcal{Z}^{\leq-1}$,
\[
\ell(z_{-\infty}^{-1}) := \inf\{k\geq1:z_{-k}=1\},
\]
with the convention that $\ell(z_{-\infty}^{-1})=+\infty$ whenever $z_{-k}=0$ for all $k\geq1$.

In this case, the value of $g(1\mid z_{-\infty}^{-1})$ depends on the distance to the most recent occurrence of the symbol $1$ in the past sequence
$z_{-\infty}^{-1}$. Thus, if $\boldsymbol{Z}$ is a chain compatible with $g$, its transition probabilities depend on a past whose length varies according to the location of the most recent occurrence of $1$. In particular,
\[
g(1\mid 0_{-\infty}^{-1})=q_{\infty},
\]
corresponding to a past in which no $1$ is observed.

If no $k\in\mathbb{N}^{*}$ exists such that
\[
q_{k+1}=q_{k+2}=\cdots=q_{\infty},
\]
then the transition rule is not determined by any fixed finite portion of
the past, and hence the model has unbounded memory. Notice that, except for
the past $0_{-\infty}^{-1}$, the relevant portion of the past is finite, but
its length is not uniformly bounded. Such processes are called \emph{chains of variable-length memory} \citep{galves2008stochastic,gallo2011chains}.

For this model, existence and uniqueness of a stationary probability
measure depend on the parameters $(q_n)_{n\in\overline{\mathbb{N}}^{*}}$.
Assume first that $q_{\infty}>0$. By \cite{cenac2012}, a stationary
probability measure compatible with $g$ exists if and only if
\[
\sum_{n\geq0} c_n<\infty,
\qquad
c_0:=1,
\qquad
c_n:=\prod_{i=1}^{n}(1-q_i),\quad n\geq1.
\]
Equivalently,
\[
1+\sum_{n\geq1}\prod_{i=1}^{n}(1-q_i)<\infty.
\]
Whenever it exists, this stationary probability measure is unique.

On the other hand, if $q_{\infty}=0$, the probability measure concentrated on the all-zero
sequence is always stationary. Moreover, \cite{cenac2012} show that if
\[
\sum_{n\geq0}c_n=\infty,
\]
this is the unique stationary probability measure, whereas if
\[
\sum_{n\geq0}c_n<\infty,
\]
the model admits a one-parameter family of stationary probability measures.
\end{example}

The existence of a $g$-measure is essential in our framework, since the information-theoretic rates considered in this work are well defined under stationarity and are precisely the population quantities we aim to estimate. Moreover, together with the additional probabilistic assumptions adopted throughout this work, it provides the setting in which our plug-in estimator converges to the population parameter of interest.  Notice, however, that uniqueness is not necessary for defining these information-theoretic quantities for a fixed stationary chain; rather, uniqueness ensures that the $g$-function determines a unique stationary probability law.

For general $g$-functions, however, the existence and uniqueness of a compatible $g$-measure are not automatic. Classical sufficient conditions are often formulated in terms of continuity properties of the $g$-function {and the rate at which its dependence on the remote past vanishes \citep{doeblin1937chaines,harris1955chains,keane1972strongly, johansson2003square}. For a past $z_{-\infty}^{-1}\in\mathcal{Z}^{\leq-1}$ and $k\in\mathbb{N}^{*}$, define the \emph{local variation} of $g$ at $z_{-\infty}^{-1}$ by \[ \operatorname{var}_k \bigl(g;z_{-\infty}^{-1}\bigr) := \sup_{\substack{ w_{-\infty}^{-1}\in\mathcal{Z}^{\leq-1}\\ w_{-k}^{-1}=z_{-k}^{-1} }} \sum_{a\in\mathcal{Z}} \left| g(a\mid w_{-\infty}^{-1}) - g(a\mid z_{-\infty}^{-1}) \right|. \] We say that $g$ is \emph{continuous} at $z_{-\infty}^{-1}$ if \[ \lim_{k\to\infty} \operatorname{var}_k \bigl(g;z_{-\infty}^{-1}\bigr) =0. \] The $g$-function is said to be \emph{continuous} if it is continuous at every $z_{-\infty}^{-1}\in\mathcal{Z}^{\leq-1}$. 

We also define the \emph{global variation} of order $k$ by  
\[ 
\operatorname{var}_k(g) := \sup_{\substack{ z_{-\infty}^{-1},\,w_{-\infty}^{-1} \in\mathcal{Z}^{\leq-1}\\ z_{-k}^{-1}=w_{-k}^{-1} }} \sum_{a\in\mathcal{Z}} \left| g(a\mid z_{-\infty}^{-1}) - g(a\mid w_{-\infty}^{-1}) \right|, \qquad k\in\mathbb{N}^{*}. 
\]
Since $\mathcal{Z}$ is finite and $\mathcal{Z}^{\leq-1}$ is compact under the product topology, continuity of $g$ is equivalent to \[ \lim_{k\to\infty}\operatorname{var}_k(g)=0. \] The sequence $\{\operatorname{var}_k(g)\}_{k\geq1}$ describes how rapidly the influence of the remote past vanishes and is commonly referred to as the \emph{variation rate} (or \emph{continuity rate}) of $g$. 

When there exists more than one stationary probability measure compatible with $g$, we say that a \emph{phase transition} occurs; when there exists only one, we speak of \emph{uniqueness}. The decay of the variations of $g$ plays an important role in classical uniqueness criteria for $g$-measures.

More recently, existence and uniqueness results have also been established beyond the classical continuity regime, allowing the $g$-function to be discontinuous under suitable structural assumptions on its set of discontinuity points \citep{gallo2013non, ferreira2020non}. The framework adopted in this paper belongs to this line of work. In particular, we do not assume global continuity of $g$. Instead, the standing assumptions introduced below combine structural and probabilistic conditions, ensuring that the stationary regime under consideration is well-defined and suitable for the concentration analysis developed later.

\paragraph{Assumption $\mathrm{(I)}$} \label{cond_i}
We assume that the $g$-function $g$ is \emph{uniformly non-null}, namely,
there exists $\epsilon_0>0$ such that
\[
g(a\mid z_{-\infty}^{-1})\geq \epsilon_0, \qquad \forall\, z_{\infty}^{-1}\in\mathcal{Z}^{\leq-1}, \quad \forall\, a\in\mathcal{Z}.
\]
In other words, uniform non-nullness ensures that the transition probability
of every symbol is bounded away from zero, uniformly over all pasts.

Following the usual terminology, a $g$-function is called \emph{regular} if it is
continuous and uniformly non-null. Otherwise, it is called \emph{non-regular}.
Thus, under Assumption~$\mathrm{(I)}$, non-regularity in our framework
arises from the possible discontinuity of the $g$-function.

\paragraph{Assumption $\mathrm{(II)}$} \label{cond_ii}

We assume that the $g$-function $g$ admits a \emph{coalescent coupling-from-the-past} (CFTP) construction with an integrable backward coalescence time. The coupling-from-the-past method was introduced by \citet{propp1996exact}. Perfect simulation constructions were subsequently developed for stochastic chains with infinite memory by \citet{comets2002processes}, and later extended to broader classes of infinite-memory kernels; see, among others, \citet{desantis2012backward}, \citet{gallo2011chains}, and \citet{gallo2013perfect}.

More precisely, following the random mapping representation used in perfect
simulation constructions for chains with infinite memory \citep{comets2002processes,gallo2013perfect}, we assume that there exists an
i.i.d. driving sequence
\[
\boldsymbol{U}:=\{U_t:t\in\mathbb{Z}\}, \qquad U_t\sim\operatorname{Uniform}[0,1],
\]
and a measurable update function
\[
F:[0,1]\times\mathcal{Z}^{\leq-1}\longrightarrow\mathcal{Z}
\]
representing the transition kernel $g$, in the sense that
\[
\lambda\left(
\left\{u\in[0,1]:
F(u,z_{-\infty}^{-1})=a
\right\}
\right)
=
g(a\mid z_{-\infty}^{-1}),
\]
for every $a\in\mathcal{Z}$ and every
$z_{-\infty}^{-1}\in\mathcal{Z}^{\leq-1}$, where $\lambda$ denotes Lebesgue measure on $[0,1]$.

Starting from an arbitrary remote past and recursively applying the update function using the driving variables, let the \emph{backward coalescence time} $\theta$ denote the amount of time one has to go backward until the value produced at time $0$ is the same for all possible remote initial
pasts. This backward-coalescence formulation is standard in perfect simulation for chains with infinite memory; see, in particular, \citet{desantis2012backward}. We assume that
\[
\mathbb{P}(\theta<\infty)=1 \qquad\text{and}\qquad \mathbb{E}[\theta]<\infty.
\]

Thus, once the construction is started sufficiently far in the past, the value produced at time $0$ no longer depends on the remote initial condition. Equivalently, the present state can be determined exactly from an almost surely finite random portion of the i.i.d.\ driving sequence.
Such coalescent constructions provide perfect simulation and, under the conditions ensuring almost sure termination, yield existence and uniqueness of the stationary probability measure compatible with the kernel; see, for instance, \citet{gallo2011chains} and \citet{gallo2013perfect} for infinite-memory
settings, including classes beyond global continuity.

For regular attractive kernels, \citet{gallo2014attractive} further establish a close connection between uniqueness of the stationary compatible chain and the existence of a coupling-from-the-past perfect simulation scheme. We emphasize that, in the present work, Assumption~$\mathrm{(II)}$ is imposed directly and does not require global continuity of the $g$-function.

The stronger integrability condition
\[
\mathbb{E}[\theta]<\infty
\]
will be used quantitatively in the concentration bounds derived below.

\paragraph{Assumption $\mathrm{(III)}$} \label{cond_iii}

We also assume that the stationary process $\boldsymbol{Z}$ compatible with the $g$-function $g$ is \emph{absolutely regular} (equivalently, \emph{$\beta$-mixing}) with an exponential rate of decay. For each $k\in\mathbb{Z}$, we define
\[
\mathcal{F}_{\boldsymbol{Z}}^{\leq k} := \sigma(Z_j:j\leq k) \qquad\text{and}\qquad \mathcal{F}_{\boldsymbol{Z}}^{\geq k} := \sigma(Z_j:j\geq k).
\]
Thus, for $n\geq1$, $\mathcal{F}_{\boldsymbol{Z}}^{\leq k}$ contains the information carried by the process up to time $k$, whereas $\mathcal{F}_{\boldsymbol{Z}}^{\geq k+n}$ contains the information carried by the process from time $k+n$ onward.

The $\beta$-mixing coefficients of $\boldsymbol{Z}$ are defined by
\[
\beta(n) := \sup_{k\in\mathbb{Z}} \beta\left(\mathcal{F}_{\boldsymbol{Z}}^{\leq k},
\mathcal{F}_{\boldsymbol{Z}}^{\geq k+n} \right), \qquad n\geq1,
\]
where
\[
\beta\left(\mathcal{F}_{\boldsymbol{Z}}^{\leq k},\mathcal{F}_{\boldsymbol{Z}}^{\geq k+n}\right):=\sup\left\{\frac{1}{2}\sum_{i=1}^{I}\sum_{j=1}^{J}\left|\mathbb{P}(A_i\cap B_j)
-\mathbb{P}(A_i)\mathbb{P}(B_j)\right|\right\},
\]
with the supremum being taken over all finite partitions
$\{A_1,\ldots,A_I\}$ and $\{B_1,\ldots,B_J\}$ of $\Omega$ such that
$A_i\in\mathcal{F}_{\boldsymbol{Z}}^{\leq k}, i=1,\ldots,I,$
and $B_j\in\mathcal{F}_{\boldsymbol{Z}}^{\geq k+n},j=1,\ldots,J.$

We say that $\boldsymbol{Z}$ is \emph{absolutely regular} if
\[
\beta(n)\longrightarrow0,
\qquad\text{as }n\to\infty.
\]
In addition, we assume that there exist positive constants $\delta$ and $\gamma$ such that
\[
\beta(n)\leq\gamma e^{-\delta n},
\qquad n\geq1.
\]

Intuitively, this condition means that the dependence between the distant
past and the distant future decreases exponentially fast as the gap between
them increases.

A non-Markovian $g$-function satisfying the Assumptions (I) - (III) is given by the following example. This example will also be used throughout the paper to illustrate the assumptions and the quantitative results obtained below.

\begin{example} 
\label{ex:alternating_case}

Consider the $g$-function defined in Example~\ref{ex:non-Markovian-g-function}. Let \[ q_{2n-1}=\frac{2}{5}, \qquad q_{2n}=\frac{1}{3}, \qquad n\geq1, \] and set \[ q_{\infty}=\frac{1}{3}. \] Thus, \[ g(1\mid z_{-\infty}^{-1}) := q_{\ell(z_{-\infty}^{-1})}, \qquad g(0\mid z_{-\infty}^{-1}) := 1-g(1\mid z_{-\infty}^{-1}). \] 

Note that $g$ is not Markov of any finite order. Indeed, for every $m\geq1$, one can choose two pasts having the same last $m$ symbols, all equal to $0$, but such that the most recent occurrence of $1$ is at distance $m+1$ in one past and at distance $m+2$ in the other. Since $m+1$ and $m+2$ have different parity, one of the corresponding transition probabilities is $2/5$ and the other is $1/3$. Hence no fixed finite portion of the past determines the transition rule. 

Moreover, $g$ is discontinuous at the all-zero past $0_{-\infty}^{-1}$. Indeed, for each $k\geq0$, consider $w^{(k)}\in\mathcal{Z}^{\leq-1}$ defined by
\[
w^{(k)}_i
=
\begin{cases}
1, & i=-(k+1),\\
0, & i\neq -(k+1),
\end{cases}
\qquad i\leq-1.
\]
Then
\[
w^{(k)} \stackrel{k \rightarrow +\infty}{\longrightarrow} 0_{-\infty}^{-1}
\]
in the product topology.
Consequently, since $\ell\bigl(w^{(k)}\bigr)=k+1$, we have
\[
g\bigl(1\mid w^{(k)}\bigr)=q_{k+1}.
\]
The sequence $(q_{k+1})_{k\geq0}$ oscillates between $2/5$ and $1/3$ and therefore does not converge to
\[
g(1\mid 0_{-\infty}^{-1})=q_\infty=\frac13.
\]
Hence, $g$ is discontinuous at $0_{-\infty}^{-1}$.

Despite this discontinuity, the $g$-function admits a unique compatible stationary probability measure. Indeed, since
\[
q_{\infty}=\frac13>0,
\]
the existence and uniqueness criterion for the infinite comb given by \citet{cenac2012} applies. Recall that, in the notation of Example~\ref{ex:non-Markovian-g-function}, one sets
\[
c_0:=1,\qquad c_n:=\prod_{i=1}^{n}(1-q_i), \qquad n\geq1.
\]
For the present choice of the sequence $(q_n)$, we have $1-q_{2j-1}=\frac35$ and $1-q_{2j}=\frac23,$ and therefore
\[
c_{2m} = \left(\frac25\right)^m \quad \hbox{ and } \quad
c_{2m+1} = \frac35\left(\frac25\right)^m,
\]
for every $m\geq0.$ Consequently,
\[
\sum_{n\geq0}c_n
=
\sum_{m\geq0}\left(\frac25\right)^m
+
\frac35\sum_{m\geq0}\left(\frac25\right)^m
=
\frac83
<\infty.
\]
Hence, by \citet{cenac2012}, there exists a unique stationary probability measure $\mu$ compatible with $g$. In particular, $\mu$ is the unique $g$-measure associated with this transition kernel.

Moreover, Assumption~$\mathrm{(I)}$ holds with $\epsilon_0=\frac13$. Indeed,  \[ g(a\mid z_{-\infty}^{-1}) \geq \frac13, \qquad \forall\,z_{-\infty}^{-1}\in\{0,1\}^{\leq-1}, \quad \forall\,a\in\{0,1\}. \] 

 We next verify Assumption~$\mathrm{(II)}$. Let $\boldsymbol{U}=\{U_t:t\in\mathbb{Z}\}$ be an i.i.d.\ sequence of $\operatorname{Uniform}[0,1]$ random variables and consider the update rule \[ Z_t = \mathbf{1} \left\{ U_t \leq q_{\ell(Z_{-\infty}^{t-1})} \right\}. \] Since $q_{\ell(z_{-\infty}^{-1})}\geq\frac13$ for every past, whenever $U_t\leq1/3$ we necessarily have $Z_t=1$, independently of the remote past. The occurrence of such a forced $1$ regenerates the chain, since from that time onward the distance to the most recent occurrence of $1$ is completely determined. This is precisely the type of regenerative CFTP mechanism described for $\epsilon$-regular symbols in \citet{gallo2011chains}.  Let \[ \tau := \inf\left\{ m\in\mathbb{N}:U_{-m}\leq\frac13 \right\}. \] Then the backward coalescence time $\theta$ is stochastically bounded by $\tau$. Since $\tau$ has a geometric tail with parameter $1/3$, \[ \mathbb{P}(\theta\geq m) \leq \mathbb{P}(\tau\geq m) = \left(\frac23\right)^m, \qquad m\geq0. \] Consequently, \[ \mathbb{P}(\theta<\infty)=1 \qquad\text{and}\qquad \mathbb{E}[\theta]<\infty, \] and Assumption~$\mathrm{(II)}$ is satisfied. 

Finally, we verify Assumption~$\mathrm{(III)}$. Fix $k\in\mathbb{Z}$ and $n\geq1$, and write
$\mathcal{F}_{\boldsymbol{Z}}^{\leq k} := \sigma(Z_j:j\leq k)$ and  $\mathcal{F}_{\boldsymbol{Z}}^{\geq k+n} := \sigma(Z_j:j\geq k+n)$. Since $\mathcal{Z}^{\geq k+n}$ is a standard Borel space, the conditional-probability representation of the coefficient of absolute
regularity \citep[Proposition~3.22]{bradley2007introduction} yields
\[\beta\left( \mathcal{F}_{\boldsymbol{Z}}^{\leq k}, \mathcal{F}_{\boldsymbol{Z}}^{\geq k+n} \right)
= \mathbb{E}\left[ \left\|\mathcal{L}\left( Z_{k+n}^{\infty} \,\middle|\, \mathcal{F}_{\boldsymbol{Z}}^{\leq k} \right) - \mathcal{L}\left(Z_{k+n}^{\infty}\right) \right\|_{\mathrm{TV}} \right],
\]
where $\mathcal{L}(X)$ denotes the law of a random element $X$, $\mathcal{L}(X\mid\mathcal{G})$ denotes a regular conditional law of $X$ given a $\sigma$-algebra $\mathcal{G}$, and $\|\cdot\|_{\mathrm{TV}}$ denotes the total variation distance.

To bound this quantity, consider a second copy $\boldsymbol{Z}'$ of the chain whose past up to time $k$ has the stationary distribution and is independent of $\mathcal{F}_{\boldsymbol{Z}}^{\leq k}$. From time $k+1$ onward, construct $\boldsymbol{Z}$ and $\boldsymbol{Z}'$ using the same driving variables $U_{k+1},U_{k+2},\ldots$. Conditionally on $\mathcal{F}_{\boldsymbol{Z}}^{\leq k}$, $Z_{k+n}^{\prime\,\infty}$ has the same distribution as $Z_{k+n}^{\infty}$ under its stationary law. Hence, by the coupling inequality \citep[Chapter~I, Section~2]{lindvall2002lectures},
\[
\left\|\mathcal{L}\left(Z_{k+n}^{\infty}\,\middle|\,\mathcal{F}_{\boldsymbol{Z}}^{\leq k}\right)
- \mathcal{L}\left(Z_{k+n}^{\infty}\right)\right\|_{\mathrm{TV}} \leq\mathbb{P}\left(Z_{k+n}^{\infty}
\neq Z_{k+n}^{\prime\,\infty}\,\middle|\,\mathcal{F}_{\boldsymbol{Z}}^{\leq k}\right) \quad a.s.
\]
If $U_{k+r}\leq1/3$ for some $r\in\{1,\ldots,n\}$, then $Z_{k+r}=Z'_{k+r}=1$
regardless of the two pasts. Since the transition rule depends only on the distance to the most recent occurrence of $1$, and both chains use the same driving variables thereafter, they coincide at all subsequent times. Consequently,
\[
\left\{Z_{k+n}^{\infty}\neq Z_{k+n}^{\prime\,\infty}\right\} \subseteq \bigcap_{r=1}^{n} \left\{ U_{k+r}>\frac13\right\}.
\]
Since the driving variables are i.i.d.\ and independent of $\mathcal{F}_{\boldsymbol{Z}}^{\leq k}$,
\[
\beta\left(
\mathcal{F}_{\boldsymbol{Z}}^{\leq k},
\mathcal{F}_{\boldsymbol{Z}}^{\geq k+n}
\right) \leq \mathbb{E}\left[\mathbb{P}\left(Z_{k+n}^{\infty}\neq Z_{k+n}^{\prime\,\infty} \,\middle|\, \mathcal{F}_{\boldsymbol{Z}}^{\leq k} \right)\right] \leq \left(\frac23\right)^n = e^{-n\log(3/2)}.
\]
Hence, Assumption~$\mathrm{(III)}$ holds, for instance, with
\[
\gamma=1,
\qquad
\delta=\log\left(\frac32\right).
\]

Therefore, Assumptions~$\mathrm{(I)}$, $\mathrm{(II)}$, and $\mathrm{(III)}$ are all satisfied.

\end{example}

\begin{remark}
\label{rmk:alternating_case}
The particular values $2/5$ and $1/3$ in Example~\ref{ex:alternating_case} are chosen only for concreteness. Indeed, the same construction can be carried out for any $0<b<a<1$ by setting
\[
q_{2n-1}=a,\qquad q_{2n}=b,\qquad n\geq1,
\]
and $q_\infty=b.$ 

The resulting $g$-function is of infinite order and is discontinuous at the all-zero past $0_{-\infty}^{-1}$.  The stationary regime is also well defined throughout this family. Indeed, \[ \sum_{n\geq0} c_n = \frac{2-a}{1-(1-a)(1-b)} <\infty. \] Since $q_\infty=b>0$, the criterion of \citet{cenac2012} implies that there exists a unique stationary probability measure compatible with $g$.

Moreover, Assumption~$\mathrm{(I)}$ holds with $\epsilon_0=\min\{b,1-a\}>0.$ Since $g(1\mid z_{-\infty}^{-1})\geq b$ for every past, the event $\{U_t\leq b\}$ forces $Z_t=1$
independently of the remote past. Consequently, the same regenerative CFTP
argument yields a backward coalescence time with geometric tail,
\[
\mathbb{P}(\theta\geq m)\leq(1-b)^m,
\qquad m\geq0,
\]
and hence Assumption~$\mathrm{(II)}$ is satisfied. The corresponding coupling argument also gives
\[
\beta(n)\leq(1-b)^n = e^{-n[-\log(1-b)]}, \qquad n\geq1,
\]
so that Assumption~$\mathrm{(III)}$ holds, for instance, with
\[\gamma=1 \quad \hbox{and} \quad \delta=-\log(1-b).
\]

Thus, Example~\ref{ex:alternating_case} is a particular member of a two-parameter family of discontinuous chains with unbounded memory satisfying Assumptions~$\mathrm{(I)}$--$\mathrm{(III)}$.
\end{remark}

\subsection{Entropy-based quantities}
\label{sec:entropy_quantities}

In this section, we introduce the basic entropy quantities associated with the stationary process under study, namely joint entropy, conditional entropy, and the entropy rate. These notions are defined with respect to a stationary probability measure $\mu$ compatible with the $g$-function $g$. These quantities will serve as the building blocks for the information-theoretic functionals considered later in the paper.

From this point on, we let $\boldsymbol{Z}$ be a stationary stochastic chain taking values in a finite alphabet $\mathcal{Z}$, with law $\mu$ compatible with the $g$-function $g$. Under this setting, we now introduce the basic entropy quantities associated with $\boldsymbol{Z}$.

Fix a positive integer $k$. The \emph{joint entropy} of the random vector $Z_{-k+1}^0$, also called the \emph{$k$-block entropy} of $\boldsymbol{Z}$, is defined by \[ H(Z_{-k+1}^0) := -\sum_{z_{-k+1}^0 \in \mathcal{Z}^k} \mu([z_{-k+1}^0]) \log \mu([z_{-k+1}^0]), \] where, throughout the paper, $\log$ denotes the natural logarithm. We adopt the usual convention $0\log 0:=0$. Thus, the entropy of a random vector quantifies the average degree of uncertainty associated with its possible outcomes.

Next, for $k\geq 2$, the entropy of $Z_0$ conditioned on the $(k-1)$-block $Z_{-k+1}^{-1}$, also called the \emph{conditional $k$-block entropy} of $\boldsymbol{Z}$, is defined as \[ H\!\left(Z_0 \mid Z_{-k+1}^{-1}\right) := -\sum_{z_{-k+1}^{0}\in\mathcal{Z}^{k}} \mu\!\left([z_{-k+1}^{0}]\right) \log\!\left( \frac{\mu\!\left([z_{-k+1}^{0}]\right)} {\mu\!\left([z_{-k+1}^{-1}]\right)} \right). \] Since $[z_{-k+1}^{0}]\subseteq[z_{-k+1}^{-1}]$, the denominator in the ratio above can vanish only when the numerator also vanishes. We therefore adopt the usual information-theoretic convention $ 0\log\frac{0}{q}:=0,\qquad q\geq0$, which, in particular, gives $0\log(0/0):=0$; see, e.g., \citep{cover1991information}. For $k=1$, we use the convention $H(Z_0\mid Z_{-k+1}^{-1})=H(Z_0)$, that is, conditioning on the empty past leaves the entropy unchanged. Equivalently, \[ H\!\left(Z_0 \mid Z_{-k+1}^{-1}\right) = H\!\left(Z_{-k+1}^0\right)-H\!\left(Z_{-k+1}^{-1}\right). \] Therefore, conditional entropy quantifies the average uncertainty that remains about the present symbol $Z_0$ after the recent past $Z_{-k+1}^{-1}$ has been observed.

The \emph{conditional entropy given the whole past} is understood as the $\mu$-average of the entropy of the conditional distribution of $Z_0$ given $Z_{-\infty}^{-1}$, namely 
{\footnotesize{
\[
H\!\left(Z_0 \mid Z_{-\infty}^{-1}\right) := -\int_{\mathcal{Z}^{\mathbb Z}}\sum_{z_0\in\mathcal{Z}}\mathbb{E}_{\mu}\!\left[\mathbf{1}_{[z_0]}\mid \mathcal{F}^{\leq -1}\right](z)\log\mathbb{E}_{\mu}\!\left[\mathbf{1}_{[z_0]}\mid \mathcal{F}^{\leq -1}\right](z)\,d\mu(z).
\]
}}
Since $\mu$ is compatible with the $g$-function $g$,
\[\mathbb{E}_{\mu}\!\left[\mathbf{1}_{[z_0]}\mid \mathcal{F}^{\leq -1}\right](z)=g\!\left(z_0\mid z_{-\infty}^{-1}\right)\qquad \mu\text{-a.s.}
\]
Therefore,
\[
H\!\left(Z_0 \mid Z_{-\infty}^{-1}\right)=-\int_{\mathcal Z^{\mathbb Z}}\sum_{z_0\in\mathcal Z}g\!\left(z_0\mid z_{-\infty}^{-1}\right)\log g\!\left(z_0\mid z_{-\infty}^{-1}\right)\,d\mu(z).
\]

Since conditioning on additional information cannot increase entropy, the sequence $\{H(Z_0\mid Z_{-k+1}^{-1})\}_{k\geq 1}$ is non-increasing and bounded below by zero. Therefore, it converges, and we define 
\[ h(\boldsymbol{Z}) := \lim_{k\to\infty} H\!\left(Z_0 \mid Z_{-k+1}^{-1}\right). \]

Moreover, for $k\geq2$, let $\mathcal{F}_{-k+1}^{-1} := \sigma(\pi_j:-k+1\leq j\leq-1).$ Since $\mathcal F_{-k+1}^{-1}\uparrow\mathcal F^{\leq-1}$ as $k\to\infty$, L\'evy's upward theorem for conditional expectations \cite[Section~14.2]{williams1991probability} implies that, for every $z_0\in\mathcal Z$,
\[
\mathbb{E}_{\mu}\!\left[\mathbf{1}_{[z_0]}\mid\mathcal{F}_{-k+1}^{-1}\right]\stackrel{k \rightarrow+\infty}{\longrightarrow}\mathbb{E}_{\mu}\!\left[\mathbf{1}_{[z_0]}\mid\mathcal{F}^{\leq-1}\right] \qquad \mu\text{-a.s.}
\]
Since $\mathcal Z$ is finite, entropy is continuous as a function of the underlying probability distribution on the probability simplex and takes values in $[0,\log|\mathcal Z|]$. Hence, by the dominated convergence theorem \cite[Section~5.9]{williams1991probability},
\[
\lim_{k\to\infty} H\!\left(Z_0\mid Z_{-k+1}^{-1}\right) = H\!\left(Z_0\mid Z_{-\infty}^{-1}\right).
\]
Therefore,
\[
h(\boldsymbol Z) = H\!\left(Z_0\mid Z_{-\infty}^{-1}\right).
\]

Finally, by the chain rule for joint entropy \cite[Theorem~2.5.1]{cover1991information} and stationarity, $$ H\!\left(Z_{-k+1}^{0}\right) = \sum_{j=1}^{k} H\!\left(Z_0\mid Z_{-j+1}^{-1}\right).$$ Hence, $\frac{1}{k}H\!\left(Z_{-k+1}^{0}\right) = \frac{1}{k}\sum_{j=1}^{k} H\!\left(Z_0\mid Z_{-j+1}^{-1}\right)$.  By Ces\`aro's theorem \cite[Theorem~4.2.3]{cover1991information}, it follows that \[ h(\boldsymbol{Z}) = \lim_{k\to\infty} H\!\left(Z_0 \mid Z_{-k+1}^{-1}\right) = \lim_{k\to\infty} \frac{1}{k}H\!\left(Z_{-k+1}^{0}\right). \] We call $h(\boldsymbol{Z})$ the \emph{entropy rate} of the chain $\boldsymbol{Z}$ \cite[Section~4.2]{cover1991information}. Intuitively, the entropy rate measures the average amount of new uncertainty produced by the process per time step, once the entire past is taken into account.

As simple illustrations, if $\boldsymbol{Z}$ is an i.i.d. chain, then
\[
h(\boldsymbol{Z}) = \lim_{k\to\infty}\frac{1}{k}H(Z_{-k+1}^{0}) = \lim_{k\to\infty}\frac{k\,H(Z_0)}{k} = H(Z_0).
\]
On the other hand, if $\boldsymbol{Z}$ is an $m$-th order Markov chain, then
for every $k\geq m+1$,
\[
H(Z_0\mid Z_{-k+1}^{-1}) = H(Z_0\mid Z_{-m}^{-1}),
\]
and therefore
\[
h(\boldsymbol{Z}) = H(Z_0\mid Z_{-m}^{-1}).
\]
Thus, in the i.i.d. case, the entropy rate coincides with the entropy of a single symbol, whereas in the Markov case it is the average uncertainty about the present symbol after conditioning on
the last $m$ symbols of the past.

\begin{example}[Entropy rate of a non-Markovian chain]
\label{ex:entropy-non-Markovian}
Consider a stationary chain compatible with the variable-length memory $g$-function of Example~\ref{ex:non-Markovian-g-function}, in the unbounded-memory regime. Since
\[
h(\boldsymbol{Z}) = H(Z_0\mid Z_{-\infty}^{-1}),
\]
the representation obtained above and the compatibility of $\mu$ with the $g$-function $g$ yield
\[
h(\boldsymbol{Z}) = -\int_{\mathcal{Z}^{\mathbb Z}} \sum_{z_0\in\{0,1\}} g(z_0\mid z_{-\infty}^{-1}) \log g(z_0\mid z_{-\infty}^{-1}) \, d\mu(z).
\]
Hence,
\[
h(\boldsymbol{Z}) = -\int_{\mathcal{Z}^{\mathbb Z}} \left[q_{\ell(z_{-\infty}^{-1})}
\log q_{\ell(z_{-\infty}^{-1})} + \bigl(1-q_{\ell(z_{-\infty}^{-1})}\bigr) \log\bigl(1-q_{\ell(z_{-\infty}^{-1}})\bigr) \right] \, d\mu(z).
\]
Therefore, in this variable-length memory case, the entropy rate is obtained by averaging the binary entropy associated with $q_{\ell(z_{-\infty}^{-1})}$ over the stationary distribution of the
infinite past.

In particular, consider the non-Markovian chain of Example~\ref{ex:alternating_case}, for which
\[
q_{2n-1}=\frac25, \qquad q_{2n}=\frac13, \qquad n\geq1, \qquad\text{and}\qquad q_\infty=\frac13.
\]
Recall from Example~\ref{ex:alternating_case} that
\[
c_{2m}=\left(\frac25\right)^m,\qquad
c_{2m+1}=\frac35\left(\frac25\right)^m,
\qquad
\sum_{n=0}^{\infty}c_n=\frac83.
\]

To determine the stationary distribution of the memory length, define
\[
A_r:=\left\{z\in\mathcal Z^{\mathbb Z}:\ell(z_{-\infty}^{-1})=r\right\}, \qquad
\nu_r:=\mu(A_r), \qquad r\geq1.
\]
For every $r\geq1$,
\[
T^{-1}A_{r+1}=A_r\cap[0].
\]
Indeed, in order that the most recent $1$ be at distance $r+1$ one time step later, it must currently be at distance $r$ and the present symbol must be $0$. Hence, by stationarity,
\[
\nu_{r+1} = \mu(A_{r+1}) = \mu(T^{-1}A_{r+1}) = \mu(A_r\cap[0]).
\]
Since $A_r\in\mathcal F^{\leq-1}$ and $\mu$ is compatible with $g$,
\[
\mu(A_r\cap[0])
=
\int_{A_r} g(0\mid z_{-\infty}^{-1})\,d\mu(z)
=
(1-q_r)\nu_r.
\]
Therefore, for every $r \geq 1$, we have $\nu_{r+1}=(1-q_r)\nu_r$. It follows recursively that $\nu_r=\nu_1 c_{r-1}$. 

Let $A_\infty := \left\{z\in\mathcal Z^{\mathbb Z}: \ell(z_{-\infty}^{-1})=\infty\right\}.$ Since $q_\infty=1/3>0$ and $T^{-1}A_\infty=A_\infty\cap[0]$, stationarity and compatibility give 
\[
\mu(A_\infty)
=
(1-q_\infty)\mu(A_\infty),
\]
and hence $\mu(A_\infty)=0$. Since the sets $A_r$, $r\geq1$, together with
$A_\infty$ form a partition of $\mathcal Z^{\mathbb Z}$, we have
\[
1=\sum_{r=1}^{\infty}\nu_r
=\nu_1\sum_{n=0}^{\infty}c_n
\Longrightarrow
\nu_1 =
\left(\sum_{n=0}^{\infty}c_n\right)^{-1},
\]
and thus
\[
\nu_r
=
\frac{c_{r-1}}
{\displaystyle\sum_{n=0}^{\infty}c_n},
\qquad r\geq1.
\]

Define the events
\[
O
:=
\bigcup_{m=0}^{\infty} A_{2m+1},
\qquad
E
:=
\bigcup_{m=0}^{\infty} A_{2m+2},
\]
corresponding, respectively, to an odd and an even finite memory length.
Since the sets $A_r$, $r\geq1$, are pairwise disjoint, we have
\[
\begin{aligned}
\mu(O)
&=
\sum_{m=0}^{\infty}\nu_{2m+1}
=
\frac{\displaystyle\sum_{m=0}^{\infty}c_{2m}}
{\displaystyle\sum_{n=0}^{\infty}c_n}
=
\frac58,\\[1ex]
\mu(E)
&=
\sum_{m=0}^{\infty}\nu_{2m+2}
=
\frac{\displaystyle\sum_{m=0}^{\infty}c_{2m+1}}
{\displaystyle\sum_{n=0}^{\infty}c_n}
=
\frac38.
\end{aligned}
\]

Since $q_{\ell(z_{-\infty}^{-1})} = \frac25 \text{ on }O, q_{\ell(z_{-\infty}^{-1})}
= \frac13 \text{ on }E,$ and $\mu(A_\infty)=0$, it follows that
\[
h(\boldsymbol Z) = \mu(O)\,h_{\mathrm b}\!\left(\frac25\right) + \mu(E)\, h_{\mathrm b}\!\left(\frac13\right)= \frac58\,h_{\mathrm b}\!\left(\frac25\right)
+ \frac38\,h_{\mathrm b}\!\left(\frac13\right),
\]
where $h_{\mathrm b}(x) := -x\log x-(1-x)\log(1-x), x\in(0,1),$ denotes the binary entropy function.
\end{example}

\begin{remark}
The computation in Example~\ref{ex:entropy-non-Markovian} extends directly
to the two-parameter family considered in Remark~1, namely
\[
q_{2n-1}=a,\qquad q_{2n}=b,\qquad n\geq1,
\qquad\text{and}\qquad q_\infty=b,
\]
with $0<b<a<1$. Using the same notation $O$ and $E$ as in the previous
example, the same argument gives
\[
\mu(O)=\frac{1}{2-a},
\qquad
\mu(E)=\frac{1-a}{2-a}.
\]
Therefore,
\[
h(\boldsymbol Z)
=
\frac{1}{2-a}\,h_{\mathrm b}(a)
+
\frac{1-a}{2-a}\,h_{\mathrm b}(b).
\]
In particular, taking $a=2/5$ and $b=1/3$ recovers the entropy rate
computed in Example~\ref{ex:entropy-non-Markovian}.
\end{remark}

\subsection{Transfer entropy}
\label{sec:transfer_entropy}

In this section, we introduce transfer entropy, the main information-theoretic quantity considered in this work. Roughly speaking, transfer entropy measures the directed flow of information from one stochastic process to another by quantifying the reduction in uncertainty about the present state of a process due to the past of another process, beyond the information already contained in the latter's own past. As such, it provides a natural measure of directional dependence between stochastic processes and will play a central role in our framework.

Let $\boldsymbol{Z}=\{(X_t,Y_t):t\in\mathbb{Z}\}$ be a stationary bivariate chain taking values in the finite alphabet $\mathcal{Z}=\mathcal{X}\times\mathcal{Y}$, with law $\mu$. For a given
positive integer $k$, the \emph{transfer entropy} from $\boldsymbol{X}:=\{X_t:t\in\mathbb{Z}\}$ to $\boldsymbol{Y}:=\{Y_t:t\in\mathbb{Z}\}$, based on blocks of length $k$, is defined by
\[
T(X_{-k+1}^0 \to Y_{-k+1}^0) := H(Y_0\mid Y_{-k+1}^{-1}) - H(Y_0\mid X_{-k+1}^{-1},Y_{-k+1}^{-1}).
\]
By the empty-block convention introduced above,
$T(X_{-k+1}^0 \to Y_{-k+1}^0)=0$ when $k=1$.

For $k\geq2$, let $\mathcal F_{-k+1}^{-1} := \sigma\!\left(\pi_j^X,\pi_j^Y:-k+1\leq j\leq-1\right)$ and $\mathcal G_{-k+1}^{-1} := \\ \sigma\!\left(\pi_j^Y:-k+1\leq j\leq-1\right),$ where
$\pi_j^X(x,y):=x_j$ and $\pi_j^Y(x,y):=y_j.$ Thus, $\mathcal F_{-k+1}^{-1}$ represents the information contained in the joint past blocks $(X_{-k+1}^{-1},Y_{-k+1}^{-1})$, whereas $\mathcal G_{-k+1}^{-1}$ represents the information contained only in the past block $Y_{-k+1}^{-1}$. For $b\in\mathcal Y$, define
\[
C_b
:=
\left\{
z=(x,y)\in\mathcal Z^{\mathbb Z}:y_0=b
\right\},
\]
i.e., $C_b$ is the event that the present symbol of the process $\boldsymbol Y$ is equal to $b$.
Then the transfer entropy can equivalently be written as
\[
T(X_{-k+1}^0\to Y_{-k+1}^0)
=
\int_{\mathcal Z^{\mathbb Z}}
\sum_{b\in\mathcal Y}
\mathbb E_\mu
\!\left[
\mathbf 1_{C_b}
\,\middle|\,
\mathcal F_{-k+1}^{-1}
\right](z)
\log
\frac{
\mathbb E_\mu
\!\left[
\mathbf 1_{C_b}
\,\middle|\,
\mathcal F_{-k+1}^{-1}
\right](z)
}{
\mathbb E_\mu
\!\left[
\mathbf 1_{C_b}
\,\middle|\,
\mathcal G_{-k+1}^{-1}
\right](z)
}
\,d\mu(z).
\]

Since conditioning on additional information cannot increase entropy \cite[Theorem 2.6.5]{cover1991information},  $T(X_{-k+1}^0 \to Y_{-k+1}^0)\geq0.$ Moreover, $T(X_{-k+1}^0\to Y_{-k+1}^0)=0$ if and only if, for every $b\in\mathcal Y$,
\[
\mathbb E_\mu\!\left[\mathbf 1_{C_b}\,\middle|\,\mathcal F_{-k+1}^{-1}\right](z) = \mathbb E_\mu\!\left[\mathbf 1_{C_b}\,\middle|\,\mathcal G_{-k+1}^{-1}\right](z),\qquad \mu\text{-a.s.}
\]
In other words, the conditional distribution of $Y_0$ given the joint past $(X_{-k+1}^{-1},Y_{-k+1}^{-1})$ coincides with its conditional distribution given only $Y_{-k+1}^{-1}$. Thus, the past of $\boldsymbol X$ provides no additional predictive information about $Y_0$ beyond that already contained in the past of $\boldsymbol Y$. In this sense, transfer entropy quantifies the directed
predictive contribution of $\boldsymbol X$ to $\boldsymbol Y$ at block length $k$.

We next consider conditioning on the whole past. Recall that $\mathcal F^{\leq-1} = \sigma(\pi_j^X,\pi_j^Y:j\leq-1)$ represents the information contained in the joint infinite past
$(X_{-\infty}^{-1},Y_{-\infty}^{-1})$. We also define $\mathcal G^{\leq-1} := \sigma(\pi_j^Y:j\leq-1),$ which represents the information contained only in the infinite past
$Y_{-\infty}^{-1}$.

Since $\mathcal Y$ is finite, there exists a measurable one-step kernel
$g_Y$ such that, for every $b\in\mathcal Y$,
\[
\mathbb E_\mu\!\left[
\mathbf 1_{C_b}
\,\middle|\,
\mathcal G^{\leq-1}
\right](z)
=
g_Y(b\mid y_{-\infty}^{-1}),
\qquad
\mu\text{-a.s.}
\]
see, e.g., \citep[Theorem 6.3]{kallenberg2002foundations}. Thus, $g_Y$ describes the conditional distribution of $Y_0$ given the
infinite past of $\boldsymbol Y$.

If $\mu$ is compatible with the bivariate $g$-function $g$, define
\[
g_{Y\mid X}
\bigl(b\mid x_{-\infty}^{-1},y_{-\infty}^{-1}\bigr)
:=
\sum_{a\in\mathcal X}
g\bigl((a,b)\mid(x_{-\infty}^{-1},y_{-\infty}^{-1})\bigr).
\]
Then, by the compatibility of $\mu$ with $g$ established in Section \ref{sec:chains_undounded_memory},
\[
\mathbb E_\mu\!\left[
\mathbf 1_{C_b}
\,\middle|\,
\mathcal F^{\leq-1}
\right](z)
=
g_{Y\mid X}
\bigl(b\mid x_{-\infty}^{-1},y_{-\infty}^{-1}\bigr),
\qquad
\mu\text{-a.s.}
\]
Since $\mathcal G^{\leq-1}\subseteq\mathcal F^{\leq-1}$, the tower property \cite[Section~9.7(i)]{williams1991probability} yields
\[
g_Y(b\mid Y_{-\infty}^{-1}) = \mathbb E_\mu\!\left[ g_{Y\mid X}\bigl(b\mid X_{-\infty}^{-1},Y_{-\infty}^{-1}\bigr)\,\middle|\,\mathcal G^{\leq-1}\right],
\qquad \mu\text{-a.s.}
\]
Consequently, the kernel $g_Y$ is, in general, not obtained by simply
summing the bivariate kernel $g$ over the present $\mathcal X$-coordinate.
Such a summation gives the conditional distribution of $Y_0$ given the
joint infinite past, whereas the conditional distribution of $Y_0$ given
only $Y_{-\infty}^{-1}$ also requires averaging over the conditional
distribution of the hidden past $X_{-\infty}^{-1}$ given
$Y_{-\infty}^{-1}$.

The \emph{transfer entropy rate} from $\boldsymbol X$ to $\boldsymbol Y$
is defined by
\[
T(\boldsymbol X\to\boldsymbol Y) := h(\boldsymbol Y)-h(\boldsymbol Y\mid\boldsymbol X),
\]
where $h(\boldsymbol Y)$ denotes the entropy rate of $\boldsymbol Y$, and
\[
h(\boldsymbol Y\mid\boldsymbol X)
:=
\lim_{k\to\infty} H(Y_0\mid X_{-k+1}^{-1},Y_{-k+1}^{-1}).
\]
Since conditioning on additional information cannot increase entropy, the sequence
$\left\{H(Y_0\mid X_{-k+1}^{-1},Y_{-k+1}^{-1})\right\}_{k\geq1}$ is non-increasing and bounded below by zero, and hence the limit exists. We refer to $h(\boldsymbol Y\mid\boldsymbol X)$ as the \emph{conditional entropy rate} of $\boldsymbol Y$ given the past of $\boldsymbol X$.

Moreover, since $\mathcal F_{-k+1}^{-1} \uparrow \mathcal F^{\leq-1}\text{ as }k\to\infty,$ the same argument based on L\'evy's upward theorem and the dominated convergence theorem used in Section~\ref{sec:entropy_quantities} yields
\[
h(\boldsymbol Y\mid\boldsymbol X) = H(Y_0\mid X_{-\infty}^{-1},Y_{-\infty}^{-1}).
\]
Using the kernel $g_{Y\mid X}$, the latter can be written as
{\footnotesize{\[
H(Y_0\mid X_{-\infty}^{-1},Y_{-\infty}^{-1}) = -\int_{\mathcal Z^{\mathbb Z}} \sum_{b\in\mathcal Y} g_{Y\mid X}\bigl(b\mid x_{-\infty}^{-1},y_{-\infty}^{-1}\bigr) \log g_{Y\mid X} \bigl(b\mid x_{-\infty}^{-1},y_{-\infty}^{-1}\bigr) \,d\mu(z).
\]}}

Since $h(\boldsymbol Y) = \lim_{k\to\infty}H(Y_0\mid Y_{-k+1}^{-1}),$ it follows that
\[
T(\boldsymbol X\to\boldsymbol Y) = \lim_{k\to\infty} T(X_{-k+1}^0\to Y_{-k+1}^0).
\]
Thus, the transfer entropy rate is the limiting directed predictive contribution of the past of $\boldsymbol X$ to $Y_0$, beyond the information already contained in the past of $\boldsymbol Y$. Equivalently, using the kernels $g_Y$ and $g_{Y\mid X}$ introduced above,
{\footnotesize{\[
T(\boldsymbol X\to\boldsymbol Y) = \int_{\mathcal Z^{\mathbb Z}} \sum_{b\in\mathcal Y} g_{Y\mid X}
\bigl(b\mid x_{-\infty}^{-1},y_{-\infty}^{-1}\bigr) \log\frac{g_{Y\mid X}\bigl(b\mid x_{-\infty}^{-1},y_{-\infty}^{-1}\bigr)}{g_Y(b\mid y_{-\infty}^{-1})}\,d\mu(z).
\]}}

Two particular cases help illustrate this definition. If the bivariate process $\boldsymbol Z=(\boldsymbol X,\boldsymbol Y)$ is i.i.d., then $Y_0$ is independent of the joint past $(X_{-\infty}^{-1},Y_{-\infty}^{-1})$, and therefore $T(\boldsymbol X\to\boldsymbol Y)=0.$ On the other hand, if the bivariate process $\boldsymbol Z=(\boldsymbol X,\boldsymbol Y)$ is a stationary Markov chain of order $m$, and if $\boldsymbol Y$ is also Markov of order $m$, then $H(Y_0\mid X_{-\infty}^{-1},Y_{-\infty}^{-1}) = H(Y_0\mid X_{-m}^{-1},Y_{-m}^{-1})$ and $H(Y_0\mid Y_{-\infty}^{-1}) = H(Y_0\mid Y_{-m}^{-1}).$ Consequently,
\[
T(\boldsymbol X\to\boldsymbol Y) = H(Y_0\mid Y_{-m}^{-1}) - H(Y_0\mid X_{-m}^{-1},Y_{-m}^{-1}),
\]
showing that, in this finite-order Markov setting, the transfer entropy rate can be expressed in terms of finite-dimensional conditional distributions.

We now introduce a motivating example, to which we will return in Section~\ref{sec:simulation} to illustrate our main results.

\begin{example}[Transfer entropy rate of a non-Markovian chain]
\label{ex:transfer-entropy-non-Markovian}

We now construct a bivariate non-Markovian chain $\boldsymbol Z=(\boldsymbol X,\boldsymbol Y)$ for which the marginal process $\boldsymbol Y$ is governed by the $g$-function of Example~\ref{ex:alternating_case}. Let $\mathcal{Z} := \mathcal X \times \mathcal Y := \{0,1\}^2$  and recall that
\[
q_{2n-1}=\frac25,
\qquad
q_{2n}=\frac13,
\qquad n\geq1,
\qquad\text{and}\qquad
q_\infty=\frac13.
\]
Fix $\eta\in(0,1/3)$. The transition mechanism of the bivariate process $\boldsymbol Z=(\boldsymbol X,\boldsymbol Y)$ is specified by the $g$-function $g:\mathcal Z^{\leq-1}\times\mathcal Z\longrightarrow[0,1].$ For a past $(x_{-\infty}^{-1},y_{-\infty}^{-1})\in\mathcal Z^{\leq-1}$, define
\[
g\bigl((x_0,y_0)\mid(x_{-\infty}^{-1},y_{-\infty}^{-1})\bigr) := \frac12\, p_\eta\bigl(y_0\mid x_{-1},y_{-\infty}^{-1}\bigr),
\]
where
\[
p_\eta\bigl(1\mid x_{-1},y_{-\infty}^{-1}\bigr) := q_{\ell(y_{-\infty}^{-1})} +\eta(2x_{-1}-1).
\]
Since $\eta<1/3$, these quantities belong to $(0,1)$ for every past, and hence $g$ defines a probability kernel on $\mathcal Z$.

Using the notation introduced above, the conditional kernel of $Y_0$ given the joint infinite past is
\[
\begin{aligned}
g_{Y\mid X}
\bigl(1\mid x_{-\infty}^{-1},y_{-\infty}^{-1}\bigr)
&=
\sum_{x_0\in\{0,1\}}
g\bigl((x_0,1)\mid(x_{-\infty}^{-1},y_{-\infty}^{-1})\bigr)
\\
&=
q_{\ell(y_{-\infty}^{-1})}
+\eta(2x_{-1}-1).
\end{aligned}
\]
As in Example~\ref{ex:alternating_case}, the dependence on $\ell(y_{-\infty}^{-1})$ cannot be determined from any fixed finite suffix of the joint past. Hence, the bivariate transition mechanism
has unbounded memory.

Let $\mu$ be a stationary probability measure compatible with $g$; its existence and uniqueness follow from the regenerative construction described at the end of the example. For every $x_0\in\{0,1\}$,
\[
\sum_{y_0\in\{0,1\}}g\bigl((x_0,y_0)\mid(x_{-\infty}^{-1},y_{-\infty}^{-1})\bigr) = \frac12.
\]
Hence, $\boldsymbol X$ is an i.i.d. Bernoulli$(1/2)$ process. Moreover,
since
\[
g\bigl((x_0,y_0)\mid(x_{-\infty}^{-1},y_{-\infty}^{-1})\bigr) = \frac12\, g_{Y\mid X} \bigl(y_0\mid x_{-\infty}^{-1},y_{-\infty}^{-1}\bigr),
\]
$X_0$ is conditionally independent of $Y_0$ given the joint infinite past. Since $X_0$ is also independent of this past, it follows that $X_0$ is independent of $Y_{-\infty}^{0}$. By stationarity,
$X_{-1}$ is therefore independent of $Y_{-\infty}^{-1}$. Consequently, using the tower property,
\[
g_Y(1\mid Y_{-\infty}^{-1})
=
\mathbb E_\mu\!\left[
g_{Y\mid X}
\bigl(1\mid X_{-\infty}^{-1},Y_{-\infty}^{-1}\bigr)
\,\middle|\,
\mathcal G^{\leq-1}
\right],
\qquad \mu\text{-a.s.}
\]
Since $X_{-1}$ is Bernoulli$(1/2)$ and independent of $Y_{-\infty}^{-1}$, for $\mu$-almost every past
$y_{-\infty}^{-1}$,
\[
\begin{aligned}
g_Y(1\mid y_{-\infty}^{-1}) = \frac12 \left(q_{\ell(y_{-\infty}^{-1})}+\eta\right) + \frac12 \left(
q_{\ell(y_{-\infty}^{-1})}-\eta\right) = q_{\ell(y_{-\infty}^{-1})}.
\end{aligned}
\]

Thus, the marginal process $\boldsymbol Y$ is governed by exactly the same $g$-function as the chain in Example~\ref{ex:alternating_case}. In particular, $\boldsymbol Y$ is non-Markovian and, by Example~\ref{ex:entropy-non-Markovian},
\[
h(\boldsymbol Y) = \frac58\, h_{\mathrm b}\!\left(\frac25\right) + \frac38\, h_{\mathrm b}\!\left(\frac13\right).
\]

Since the marginal process $\boldsymbol Y$ is governed by the same
$g$-function as in Example~\ref{ex:alternating_case}, Example~\ref{ex:entropy-non-Markovian} shows that,
\[
q_{\ell(Y_{-\infty}^{-1})} = \frac25 \quad\text{with probability }\frac58, \qquad q_{\ell(Y_{-\infty}^{-1})} = \frac13 \quad\text{with probability }\frac38.
\]
Using the expression for $g_{Y|X}$ obtained above, we obtain 
{\footnotesize{
\[
h(\boldsymbol Y\mid\boldsymbol X) ={} \frac58 \left[\frac12\,h_{\mathrm b}\!\left(\frac25+\eta\right) + \frac12\,h_{\mathrm b}\!\left(\frac25-\eta\right)\right] + \frac38\left[\frac12\,h_{\mathrm b}\!\left(\frac13+\eta\right)+\frac12\,h_{\mathrm b}\!\left(\frac13-\eta\right)\right].
\]
}}
Therefore,
\[
\begin{aligned}
T(\boldsymbol X\to\boldsymbol Y)
={}&
\frac58
\left[
h_{\mathrm b}\!\left(\frac25\right)
-\frac12\,
h_{\mathrm b}\!\left(\frac25+\eta\right)
-\frac12\,
h_{\mathrm b}\!\left(\frac25-\eta\right)
\right]
\\
&+
\frac38
\left[
h_{\mathrm b}\!\left(\frac13\right)
-\frac12\,
h_{\mathrm b}\!\left(\frac13+\eta\right)
-\frac12\,
h_{\mathrm b}\!\left(\frac13-\eta\right)
\right].
\end{aligned}
\]
Since $h_{\mathrm b}$ is strictly concave on $(0,1)$, and since $\frac25\pm\eta,
\hbox{ and } \frac13\pm\eta$ belong to $(0,1)$, both terms in brackets above are strictly positive
whenever $\eta>0$, and therefore $T(\boldsymbol X\to\boldsymbol Y)>0.$

Finally, this example also lies within the probabilistic framework of the paper. Indeed, the kernel is uniformly non-null, with $\epsilon_0=\frac12(\frac13-\eta)>0$, and the same regenerative mechanism used in Example~\ref{ex:alternating_case} applies here, i.e., whenever the auxiliary uniform variable used to generate $Y_t$ falls below $\frac13-\eta$, the value $Y_t=1$ is forced independently of the remote past. Hence, the coupling-from-the-past construction has a geometrically decaying coalescence time, with $\mathbb P(\theta\geq m)\leq(\frac23+\eta)^m$, and the same coupling argument yields exponential $\beta$-mixing. Consequently, the bivariate chain satisfies Assumptions~(I)--(III).
\end{example}

\subsection{Plug-in estimators for entropy and transfer entropy rates}
\label{sec:plugin-estimators}

We now introduce the plug-in estimators for the entropy rate and the transfer entropy rate considered in this work. Consider positive integers $k$ and $n$ such that $k\leq n$, and let
$\boldsymbol Y=\{Y_t:t\in\mathbb Z\}$ be the canonical stationary and ergodic chain on $\mathcal Y^{\mathbb Z}$, with law $\mu_Y$, taking values in the finite alphabet $\mathcal Y$. Let $y_{-n}^{-1}\in\mathcal Y^n$ be an observed sample of $\boldsymbol Y$. Let $\widetilde y$ denote the bi-infinite periodic extension of the sample $y_{-n}^{-1}$, with period
$n$, that is, $\widetilde y_{j+qn}=y_j,$ for every $j\in\{-n,\ldots,-1\}$ and $q\in\mathbb Z$.

For every $b_{-k+1}^{0}\in\mathcal Y^k$, we define the \emph{empirical frequency} of the block $b_{-k+1}^{0}$ by
\[
\widehat\mu_{n,k} \bigl(b_{-k+1}^{0}; y_{-n}^{-1}\bigr) := \frac1n \sum_{i=-n}^{-1}\mathbf 1\left\{\widetilde y_i^{\,i+k-1}=b_{-k+1}^{0}\right\}.
\]
In other words, $\widehat\mu_{n,k}(b_{-k+1}^{0}; y_{-n}^{-1})$ is obtained by sliding the
block $b_{-k+1}^{0}$ over the periodic extension of the observed sample and computing its relative frequency. Notice that $\widehat\mu_{n,k}(\cdot; y_{-n}^{-1})$ is a probability measure on
$\mathcal Y^k$. The semicolon in the above definition emphasizes that $y_{-n}^{-1}$ is the observed sample from which the empirical distribution is constructed. When there is no risk of ambiguity, we write $\widehat\mu_{n,k}(b_{-k+1}^{0})$ for $\widehat\mu_{n,k}(b_{-k+1}^{0};y_{-n}^{-1})$.

The \emph{plug-in estimator of the $k$-block entropy} is then defined by
\[
\widehat H_n(Y_{-k+1}^{0}) := - \sum_{b_{-k+1}^{0}\in\mathcal Y^k} \widehat\mu_{n,k}(b_{-k+1}^{0}) \log \widehat\mu_{n,k}(b_{-k+1}^{0}).
\]
Accordingly, for $k\geq2$, we define the \emph{plug-in estimator of the conditional entropy} by
\begin{equation}
\widehat H_n\bigl(Y_0\mid Y_{-k+1}^{-1}\bigr) := - \sum_{b_{-k+1}^{0}\in\mathcal Y^k} \widehat\mu_{n,k}(b_{-k+1}^{0})\log\frac{\widehat\mu_{n,k}(b_{-k+1}^{0})}{\widehat\mu_{n,k-1}(b_{-k+1}^{-1})}.
\label{eq:condition_entropy}
\end{equation}
As before, we use the convention $0\log(0/q)=0$ for $q\geq0$. For $k=1$, the empty-block convention gives $\widehat H_n \bigl(Y_0\mid Y_{-k+1}^{-1}\bigr) = \widehat H_n(Y_0).$ Equivalently,
$\widehat H_n\bigl(Y_0\mid Y_{-k+1}^{-1}\bigr)=\widehat H_n(Y_{-k+1}^{0})-\widehat H_n(Y_{-k+1}^{-1}).$

We next record the consistency of this estimator when the block length $k$ is fixed. For $b_{-k+1}^{0}\in\mathcal Y^k$, define
\begin{equation}
\overline\mu_{n,k}(b_{-k+1}^{0}) := \frac1n \sum_{i=-n}^{-1} \mathbf 1 \left\{Y_i^{\,i+k-1}=b_{-k+1}^{0}\right\}.
\label{eq:auxiliary_empirical_distribution}
\end{equation}
Notice that $\overline\mu_{n,k}(b_{-k+1}^{0})$ is an auxiliary random quantity defined from the underlying process $\boldsymbol Y$, rather than solely from the observed sample. Indeed, for the last $k-1$ starting positions, the blocks $Y_i^{i+k-1}$ involve coordinates after time $-1$, which are not contained in the observed sample. This quantity is introduced because it is an ordinary ergodic average, to which Birkhoff's ergodic theorem applies directly. In this context, since $\boldsymbol Y$ is stationary and ergodic, Birkhoff's ergodic  theorem \cite[Section~I.3]{shields1996ergodic} implies that, for every fixed $k\geq1$ and every $b_{-k+1}^{0}\in\mathcal Y^k$,
\[\overline\mu_{n,k}(b_{-k+1}^{0}) \stackrel{n \rightarrow +\infty}{\longrightarrow} \mu_Y([b_{-k+1}^{0}]), \qquad \mu_Y\text{-a.s.}
\]
The empirical frequencies based on the periodic extension can differ from $\overline\mu_{n,k}$ only for the last $k-1$ starting positions. Therefore,
\[\left|\widehat\mu_{n,k}(b_{-k+1}^{0}\bigr) -\overline\mu_{n,k}(b_{-k+1}^{0})\right|\leq\frac{k-1}{n}.
\]
Hence, for every fixed $k\geq1$,
\[\widehat\mu_{n,k}\bigl(b_{-k+1}^{0}; Y_{-n}^{-1}\bigr) \stackrel{n \rightarrow +\infty}{\longrightarrow} \mu_Y([b_{-k+1}^{0}]), \qquad \mu_Y\text{-a.s.}
\]
Since the alphabet is finite and entropy is a continuous function of the
underlying probability distribution, it follows that, for every fixed $k$,
\[
\widehat H_n\bigl(Y_0\mid Y_{-k+1}^{-1}\bigr) \stackrel{n \rightarrow +\infty}{\longrightarrow} H \bigl(Y_0\mid Y_{-k+1}^{-1}\bigr), \qquad \mu_Y\text{-a.s.}
\]
Consequently, by the characterization of the entropy rate obtained in Section~\ref{sec:entropy_quantities},
\[\lim_{k\to\infty}\lim_{n\to\infty} \widehat H_n \bigl(Y_0\mid Y_{-k+1}^{-1}\bigr) = h(\boldsymbol Y), \qquad \mu_Y\text{-a.s.}
\]

We next define the plug-in estimator of transfer entropy. Let
$\boldsymbol Z=(\boldsymbol X,\boldsymbol Y)$ be a stationary and ergodic
bivariate chain taking values in
$\mathcal Z=\mathcal X\times\mathcal Y$, with law $\mu$, and let
$(x_{-n}^{-1},y_{-n}^{-1})\in\mathcal X^n\times\mathcal Y^n$ be a
finite sample. Let $(\widetilde x,\widetilde y)$ denote its bi-infinite
periodic extension.

For every $(a_{-k+1}^{0},b_{-k+1}^{0})\in\mathcal X^k\times\mathcal Y^k$, define
\[
\widehat\mu_{n;k,k}\bigl(a_{-k+1}^{0},b_{-k+1}^{0}; x_{-n}^{-1},y_{-n}^{-1} \bigr) :=
\frac1n \sum_{i=-n}^{-1} \mathbf 1 \left\{\widetilde x_i^{\,i+k-1}=a_{-k+1}^{0},\, \widetilde y_i^{\,i+k-1}=b_{-k+1}^{0}\right\}.
\]
Thus, $\widehat\mu_{n;k,k}(\cdot,\cdot;x_{-n}^{-1},y_{-n}^{-1})$ is the empirical joint distribution of the length-$k$ blocks of $\boldsymbol X$ and $\boldsymbol Y$, obtained by sliding the two blocks
synchronously over the periodic extension of the observed sample. Throughout the paper, whenever the observed sample is clear from the context, we use the shorthand notation $\widehat\mu_{n;k,k}(a_{-k+1}^{0},b_{-k+1}^{0}) := \widehat\mu_{n;k,k} \bigl(a_{-k+1}^{0},b_{-k+1}^{0};x_{-n}^{-1},y_{-n}^{-1}\bigr)$.

For $0\leq r,s\leq k$, we denote by $\widehat\mu_{n;r,s}$ the corresponding empirical marginal distribution, namely
\[
\widehat\mu_{n;r,s}\bigl(a_{-k+1}^{-k+r},b_{-k+1}^{-k+s};x_{-n}^{-1},y_{-n}^{-1}\bigr)
:=\sum_{\substack{a_{-k+r+1}^{0}\in\mathcal X^{k-r}\\b_{-k+s+1}^{0}\in\mathcal Y^{k-s}}}
\widehat\mu_{n;k,k}\bigl(a_{-k+1}^{0},b_{-k+1}^{0};x_{-n}^{-1},y_{-n}^{-1}\bigr).
\]
with the usual empty-block convention when $r=0$ or $s=0$. We use the same shorthand convention for these empirical marginals.  In particular, the empirical conditional entropy of $Y_0$ given only its $(k-1)$-past defined in Equation \eqref{eq:condition_entropy} may be written as
\[
\widehat H_n\bigl(Y_0\mid Y_{-k+1}^{-1}\bigr) = - \sum_{b_{-k+1}^{0}\in\mathcal Y^k}\widehat\mu_{n;0,k}(b_{-k+1}^{0})\log\frac{\widehat\mu_{n;0,k}(b_{-k+1}^{0})}{\widehat\mu_{n;0,k-1}(b_{-k+1}^{-1})
}.
\]
Similarly, the \emph{plug-in estimator of the entropy of $Y_0$ given the joint
$(k-1)$-past} may be defined as
{\footnotesize{\[
\widehat H_n\bigl(Y_0\mid X_{-k+1}^{-1},Y_{-k+1}^{-1}\bigr) := -\sum_{\substack{a_{-k+1}^{-1}\in\mathcal X^{k-1}\\ b_{-k+1}^{0}\in\mathcal Y^k
}}
\widehat\mu_{n;k-1,k}
\bigl(a_{-k+1}^{-1},b_{-k+1}^{0}\bigr)
\log
\frac{
\widehat\mu_{n;k-1,k}
\bigl(a_{-k+1}^{-1},b_{-k+1}^{0}\bigr)
}{
\widehat\mu_{n;k-1,k-1}
\bigl(a_{-k+1}^{-1},b_{-k+1}^{-1}\bigr)
}.
\]}}
We then define the \emph{plug-in estimator of the transfer entropy at block length $k$} by
\[
\widehat T_n^k\bigl(X_{-k+1}^{0}\to Y_{-k+1}^{0}\bigr) := \widehat H_n \bigl(Y_0\mid Y_{-k+1}^{-1}\bigr) -
\widehat H_n\bigl(Y_0\mid X_{-k+1}^{-1},Y_{-k+1}^{-1}\bigr).
\]
By the empty-block convention, $\widehat T_n^1(X_0\to Y_0)=0$.

For every fixed $k$, the same argument used above, now applied to the indicator functions of joint cylinder events, shows that all the empirical frequencies appearing in the definition of $\widehat T_n^k$ converge $\mu$-almost surely to their corresponding stationary block probabilities. Indeed, the periodic and non-periodic empirical frequencies can differ only at the last $k-1$ starting positions, and hence their difference is bounded in absolute value by $(k-1)/n$. Therefore, since the alphabets are finite,
\[
\widehat T_n^k\bigl(X_{-k+1}^{0}\to Y_{-k+1}^{0}\bigr)
\stackrel{n \rightarrow +\infty}{\longrightarrow}
T\bigl(X_{-k+1}^{0}\to Y_{-k+1}^{0}\bigr), \qquad \mu\text{-a.s.}
\]
Combining this with the characterization of the transfer entropy rate obtained above, we obtain
\[
\lim_{k\to\infty} \lim_{n\to\infty} \widehat T_n^k\bigl(X_{-k+1}^{0}\to Y_{-k+1}^{0}\bigr) =
T(\boldsymbol X\to\boldsymbol Y), \qquad \mu\text{-a.s.}
\]

In order to estimate the corresponding rates, we shall allow the block length to increase with the sample size. Thus, we consider a sequence of positive integers $\{k(n)\}_{n\geq1}$ such that $k(n) \rightarrow \infty$ and, for all sufficiently large $n$, $k(n)\leq C\log n,$
for some constant $C>0$. The precise range of $C$ required for the
concentration results depends on the estimator under consideration and
will be specified in Section~\ref{sec:main_results}.

\section{Main results}
\label{sec:main_results}

We first consider the plug-in estimator of the entropy rate. Our first result provides a non-asymptotic concentration inequality around its expectation under Assumption~{\rm (II)} alone (Theorem \ref{thm:concentration_entropy_mean}). We then specialize this bound to logarithmically growing block lengths and, under the additional Assumptions~{\rm (I)} and~{\rm (III)}, control the bias (Remark \ref{rem:entropy_logarithmic_blocks}) in order to obtain concentration around the entropy rate itself (Theorem \ref{thm:concentration_entropy_rate}).

\begin{teo}
\label{thm:concentration_entropy_mean}
Let $\boldsymbol Y=\{Y_t:t\in\mathbb Z\}$ be a stationary stochastic chain on a finite alphabet $\mathcal Y$, with law $\mu_Y$, compatible with a $g$-function satisfying Assumption~{\rm (II)}. Then, for every $n\geq2$, every integer $k$ such that $2\leq k\leq n$, and every $\xi>0$,
\[
\mu_Y\left(
\left|\widehat H_n^k - E_{\mu_Y}\left[\widehat H_n^k\right]\right| > \xi\right) \leq
4\exp\left\{-\frac{2\xi^2n}{9(1+E[\theta])^2(2k-1)^2(\log n+1)^2}\right\},
\]
where, for shorthand notation, we write $\widehat H_n^k := \widehat H_n\left(Y_0\mid Y_{-k+1}^{-1}\right)$.
\end{teo}

\begin{remark}
\label{rem:entropy_logarithmic_blocks}
Theorem~\ref{thm:concentration_entropy_mean} is non-asymptotic and does not require the block length to diverge with $n$. In particular, let $\{k(n)\}_{n\geq1}$ be any sequence of positive integers such that, for all sufficiently large $n$, $2 \leq k(n)\leq C\log n,$ for some constant $C>0$. Then, for every $\xi>0$ and all sufficiently large $n$,
\[
\mu_Y\left(\left|
\widehat H_n^{k(n)} - E_{\mu_Y}\left[\widehat H_n^{k(n)}\right]\right| > \xi\right) \leq
4\exp\left\{-\frac{\xi^2n}{72C^2(1+E[\theta])^2\log^4 n}\right\},
\]
where, for shorthand notation, we write $\widehat H_n^{k(n)} := \widehat H_n\left(Y_0\mid Y_{-k(n)+1}^{-1}\right)$.
\end{remark}

\begin{teo}
\label{thm:concentration_entropy_rate}
Let $\boldsymbol Y=\{Y_t:t\in\mathbb Z\}$ be a stationary stochastic chain on a finite alphabet $\mathcal Y$, with $|\mathcal Y|\geq2$, and law $\mu_Y$, compatible with a $g$-function satisfying Assumptions~{\rm (I)}--{\rm (III)}. Let $\{k(n)\}_{n\geq1}$ be a any sequence of positive integers satisfying $k(n)\to\infty$ and, for all sufficiently large $n$, $k(n)\leq C\log n,$ $0<C<\frac{1}{\log|\mathcal Y|}.$
Then there exists a deterministic sequence $(r_n)_{n\geq1}$, with
$r_n\to0$, such that, for every $\xi>0$ and all sufficiently large $n$,
\[
\mu_Y\left(\left| \widehat H_n^{k(n)} - h(\boldsymbol Y) \right| >\xi+r_n\right)
\leq4\exp\left\{-\frac{\xi^2n}{72C^2(1+E[\theta])^2\log^4 n}\right\},
\]
where, for shorthand notation, we write $\widehat H_n^{k(n)} := \widehat H_n\left(Y_0\mid Y_{-k(n)+1}^{-1}\right)$.
\end{teo}

\begin{cor}
\label{cor:concentration_entropy_rate}
Under the assumptions of Theorem~\ref{thm:concentration_entropy_rate},
\[\widehat H_n\!\left(Y_0\mid Y_{-k(n)+1}^{-1}\right) \longrightarrow h(\boldsymbol Y)
\qquad
\mu_Y\text{-almost surely}.
\]
\end{cor}

We now turn to the plug-in estimator of the transfer entropy rate. The results follow the same structure as in the entropy case: we first establish a non-asymptotic concentration inequality around the expectation (Theorem \ref{thm:concentration_transfer_entropy_mean}) and then, after controlling the corresponding bias under Assumptions~{\rm (I)}--{\rm (III)} (Remark \ref{rmk:concentration_transfer_entropy_mean}), derive concentration around the transfer entropy rate (Theorem \ref{thm:concentration_transfer_entropy_rate}).

\begin{teo}
\label{thm:concentration_transfer_entropy_mean}
Let $\boldsymbol Z=\{(X_t,Y_t):t\in\mathbb Z\}$ be a stationary bivariate stochastic chain on a finite alphabet $\mathcal Z=\mathcal X\times\mathcal Y$, with law $\mu$, compatible with a
$g$-function satisfying Assumption~{\rm (II)}. Then, for every $n\geq2$, every integer $k$ such that $2\leq k\leq n$, and every $\xi>0$,
\[
\mu\Bigg(\left|\widehat T_n^{k} - E_\mu\!\left[\widehat T_n^{k}\right]\right|>\xi\Bigg)
\leq 4\exp\left\{-\frac{2\xi^2 n}{9(1+E[\theta])^2(4k-2)^2(\log n+1)^2}\right\},
\]
where, for shorthand notation, we write $\widehat T_n^{k} := T_n^{k}\left(X_{-k+1}^{0}\to Y_{-k+1}^{0}\right).$
\end{teo}

\begin{remark}
\label{rmk:concentration_transfer_entropy_mean}
Theorem~\ref{thm:concentration_transfer_entropy_mean} is non-asymptotic
and does not require the block length to diverge with $n$. In particular,
let $\{k(n)\}_{n\geq1}$ be any sequence of positive integers such that,
for all sufficiently large $n$, $2 \leq k(n)\leq C\log n,$ for some constant $C>0$. Then, for every $\xi>0$ and all sufficiently large $n$,
\[
\mu\Bigg(\left|\widehat T_n^{k(n)}-E_\mu\!\left[\widehat T_n^{k(n)}\right]\right|>\xi
\Bigg)\leq4\exp\left\{-\frac{\xi^2 n}{288C^2(1+E[\theta])^2\log^4 n}\right\},
\]
where, for shorthand notation, we write $\widehat T_n^{k(n)} := \widehat T_n^{k(n)}\left(X_{-k(n)+1}^{0}\to Y_{-k(n)+1}^{0}\right)$.
\end{remark}

\begin{teo}
\label{thm:concentration_transfer_entropy_rate}
Let $\boldsymbol Z=\{(X_t,Y_t):t\in\mathbb Z\}$ be a stationary bivariate stochastic chain on the finite alphabet $\mathcal Z=\mathcal X\times\mathcal Y$, with $|\mathcal Z|\geq2$, and law $\mu$, compatible with a $g$-function satisfying Assumptions~{\rm (I)}--{\rm (III)}. Let $\{k(n)\}_{n\geq1}$ be any sequence of positive integers satisfying $k(n)\to\infty$ and, for all sufficiently large $n$, $k(n)\leq C\log n,$ and $0<C<\frac{1}{\log|\mathcal Z|}.$ Then there exists a deterministic sequence $(r_n)_{n\geq1}$, with
$r_n\to0$, such that, for every $\xi>0$ and all sufficiently large $n$,
\[
\mu\left(\left|\widehat T_n^{k(n)} - T(\boldsymbol X\to\boldsymbol Y)\right| > \xi+r_n\right)
\leq 4\exp\left\{-\frac{\xi^2 n}{288C^2(1+E[\theta])^2\log^4 n}\right\},
\]
where, for shorthand notation, we write $\widehat T_n^{k(n)} := \widehat T_n^{k(n)}\left(X_{-k(n)+1}^{0}\to Y_{-k(n)+1}^{0}\right)$.
\end{teo}

\begin{cor}
\label{cor:concentration_transfer_entropy_rate}
Under the assumptions of Theorem~\ref{thm:concentration_transfer_entropy_rate},
\[
\widehat T_n^{k(n)}\left(X_{-k(n)+1}^{0}\to Y_{-k(n)+1}^{0}\right) \longrightarrow
T(\boldsymbol X\to\boldsymbol Y)
\qquad \mu\text{-almost surely}.
\]
\end{cor}

\section{Simulation}
\label{sec:simulation}

In this section, we illustrate the finite-sample behavior of the plug-in estimators considered in the previous sections. We first study the entropy-rate estimator for the non-Markovian chain of Example~\ref{ex:alternating_case}, thereby illustrating Theorem~\ref{thm:concentration_entropy_rate} and Corollary~\ref{cor:concentration_entropy_rate}. We then consider the plug-in estimator of the transfer entropy rate for the non-Markovian chain of Example \ref{ex:transfer-entropy-non-Markovian}, illustrating Theorem \ref{thm:concentration_transfer_entropy_rate} and Corollary \ref{cor:concentration_transfer_entropy_rate}.

\subsection{Entropy-rate estimation}
\label{sec:simulation_entropy_rate_estimation}

We consider the binary chain of Example~\ref{ex:alternating_case}, for which
\[
h(\boldsymbol Y) = \frac{5}{8}h_b\left(\frac25\right) + \frac{3}{8}h_b\left(\frac13\right)
\approx 0.659325.
\]

The Monte Carlo experiment consists of $R=500$ independent stationary trajectories, each of length $50\,000$. For each trajectory, we consider the sample lengths $n \in \{5000,10000,\ldots,50000\}$ and compute the plug-in estimator using the periodic empirical block distribution introduced in Section~\ref{sec:plugin-estimators}. The block length is chosen as $k(n) =\left\lfloor\frac{\log n}{2\log 2}\right\rfloor.$
Hence $k(n)\to\infty$ and $k(n)\leq C\log n,$  $C=\frac{1}{2\log 2} < \frac{1}{\log 2},$ so that the growth condition of Theorem~\ref{thm:concentration_entropy_rate} is satisfied.

Figure~\ref{fig:as_convergence_entropy} shows the Monte Carlo median of $\widehat H_n\left(Y_0\mid Y_{-k(n)+1}^{-1}\right),$
together with the empirical $5\%$--$95\%$ quantile band (grey region). The dashed horizontal line represents the exact value of $h(\boldsymbol Y)$. As the sample size increases, the empirical distribution becomes more concentrated and the median approaches the entropy rate. This behavior is consistent with the almost-sure convergence established in Corollary~\ref{cor:concentration_entropy_rate}. 
%The small local fluctuations are expected in a finite Monte Carlo experiment and are also affected by the discrete changes in the integer-valued sequence $k(n)$.
\newpage

\begin{figure}[http!]
    \centering
    \includegraphics[scale=.29]{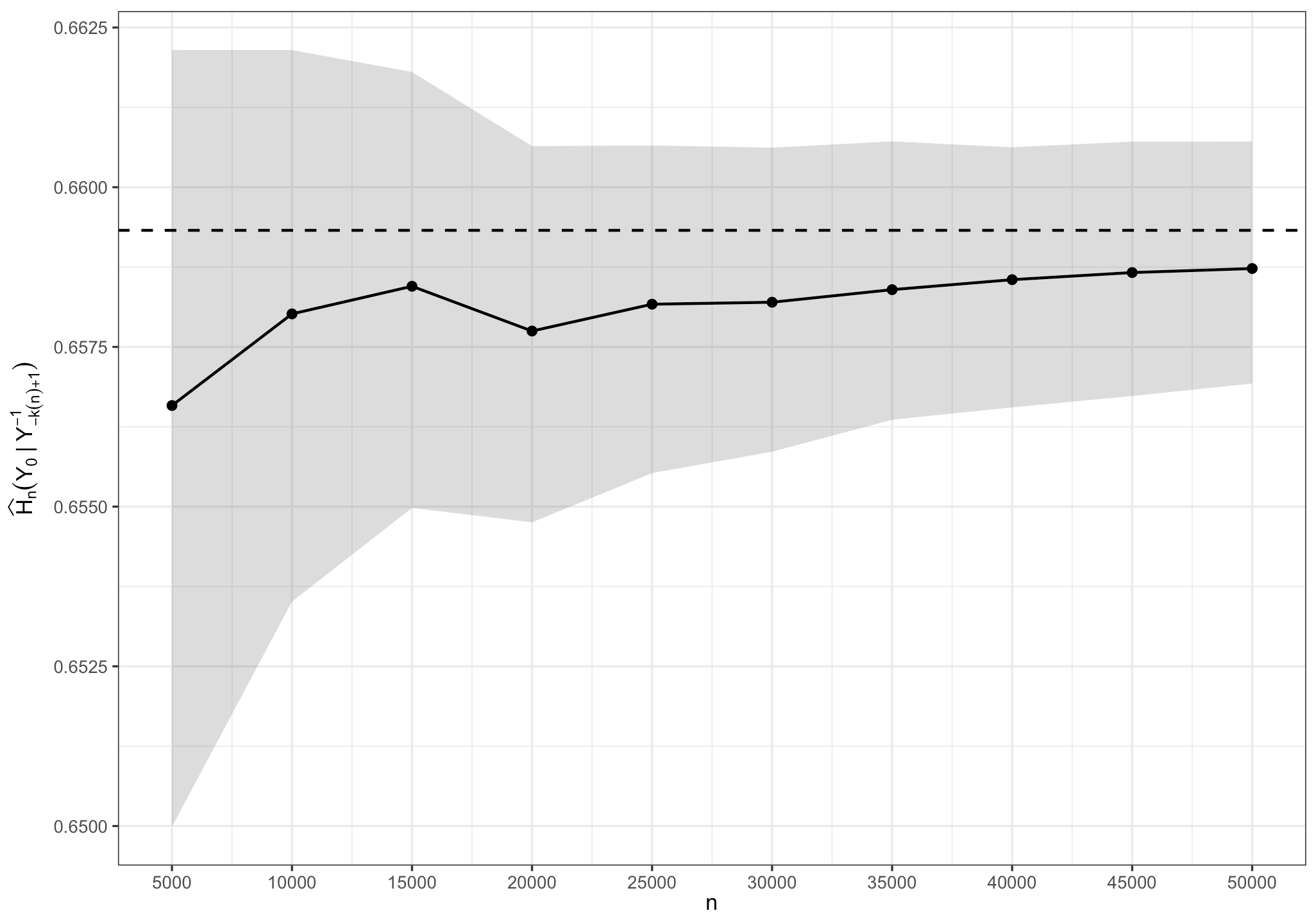}
    \caption{Convergence of the plug-in entropy-rate estimator for the chain of Example \ref{ex:alternating_case}.}
    \label{fig:as_convergence_entropy}
\end{figure}

To examine concentration more directly, Figure~\ref{fig:empirical_exceedance_entropy} displays the empirical exceedance probabilities
\[
\widehat P_n(\varepsilon) := \frac{1}{R}\sum_{r=1}^{R}\mathbf 1\left\{\left|\widehat H_n^{(r)}
-h(\boldsymbol Y)\right|>\varepsilon
\right\},
\]
for $\varepsilon\in\{0.001,0.002,0.003\}$. For each fixed deviation level, the probabilities exhibit an overall decreasing trend as $n$ increases, with a faster decrease for larger values of $\varepsilon$. This provides a finite-sample illustration of the concentration phenomenon underlying Theorem~\ref{thm:concentration_entropy_rate}. 

\begin{figure}[http!]
    \centering
    \includegraphics[scale=.29]{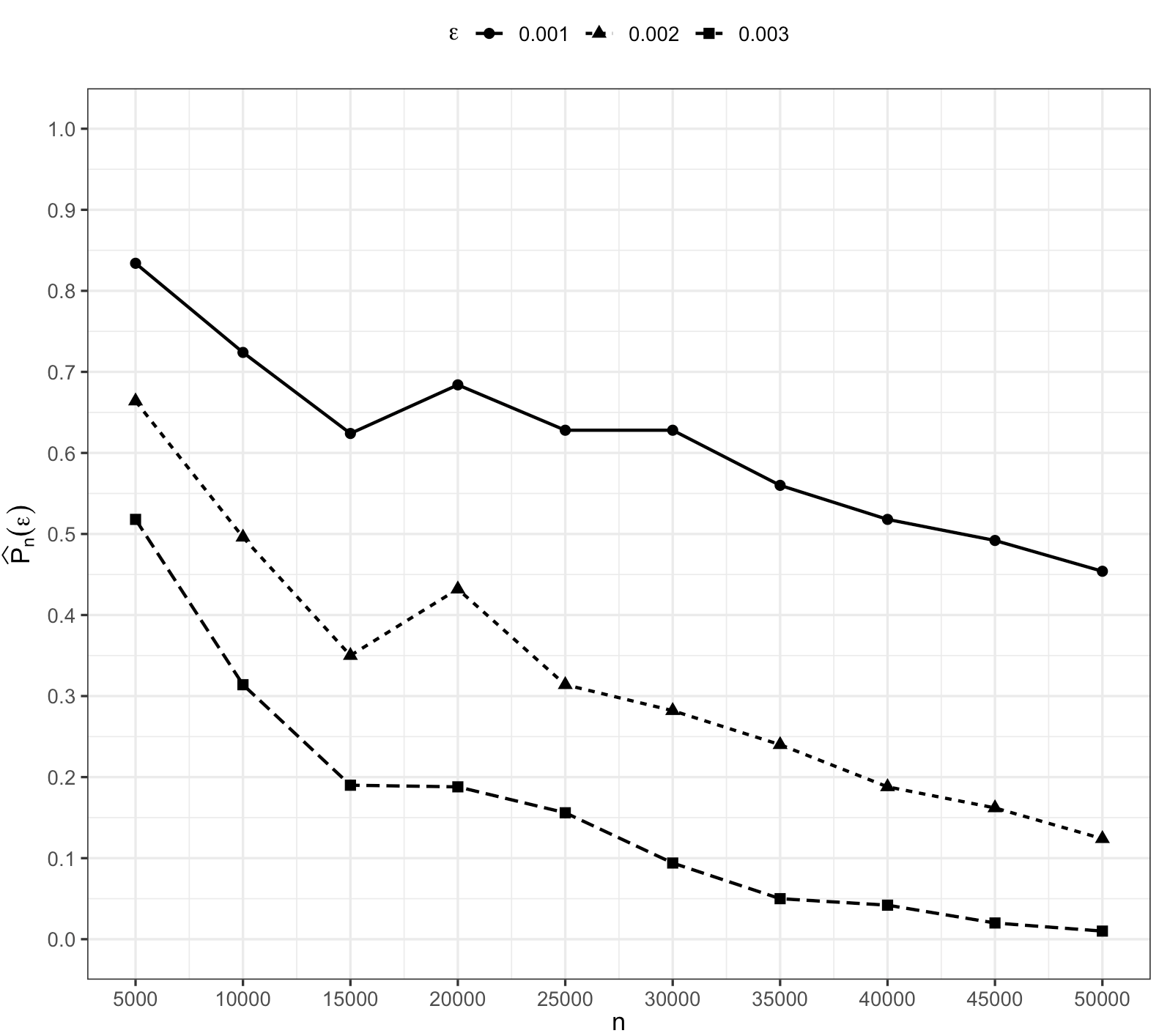}
    \caption{Empirical exceedance probabilities of the entropy-rate estimator for the chain of Example \ref{ex:alternating_case}.}
    \label{fig:empirical_exceedance_entropy}
\end{figure}

Finally, Theorem~\ref{thm:concentration_entropy_rate} gives an exponential upper bound whose relevant scale is {\footnotesize{$e^{-c\,\frac{n}{\log^4 n}.}$}} Motivated by this rate, Figure~\ref{fig:concentration_scale_entropy} plots $-\log \widehat P_n(\varepsilon)$ against $\frac{n}{\log^4 n}.$
The increasing, approximately linear trend is consistent with an exponential-type decay on the scale $n/\log^4 n$ predicted by Theorem~\ref{thm:concentration_entropy_rate}. 
%This figure should be interpreted only as a qualitative illustration of the theoretical concentration rate, since Theorem~1 provides an upper bound, rather than an exact asymptotic expression for the exceedance probabilities.

\begin{figure}[http!]
    \centering
    \includegraphics[scale=.29]{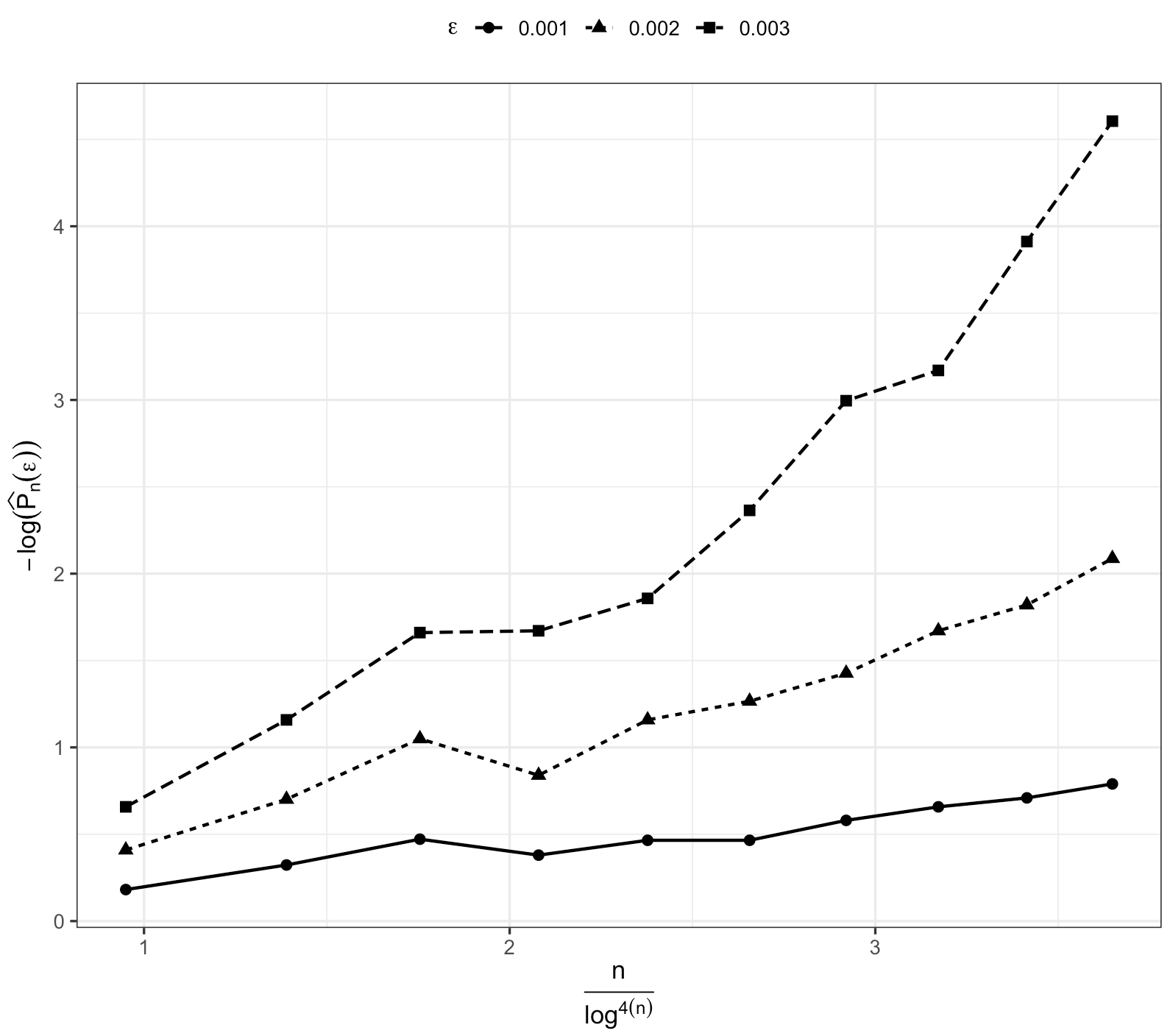}
    \caption{Exponential-type concentration scale for the chain of Example \ref{ex:alternating_case}.}
    \label{fig:concentration_scale_entropy}
\end{figure}

\subsection{Transfer-entropy-rate estimation}
\label{sec:simulation_transfer_entropy_rate_estimation}

We now consider the bivariate non-Markovian chain of Example~\ref{ex:transfer-entropy-non-Markovian}, taking $\eta=\frac{1}{5}.$ For this choice of parameter, the transfer entropy rate computed in Example~\ref{ex:transfer-entropy-non-Markovian} is
$$
\begin{aligned}
T(\boldsymbol X\to\boldsymbol Y) ={}& \frac{5}{8}\left[h_b\left(\frac{2}{5}\right)-\frac{1}{2}h_b\left(\frac{3}{5}\right)-\frac{1}{2}h_b\left(\frac{1}{5}\right)\right]\\
&+\frac{3}{8}\left[h_b\left(\frac{1}{3}\right)-\frac{1}{2}h_b\left(\frac{8}{15}\right)
-\frac{1}{2}h_b\left(\frac{2}{15}\right)\right] \approx 0.089459.
\end{aligned}
$$

As in the entropy-rate experiment, the Monte Carlo simulation consists of $R=500$ independent stationary trajectories, each of length $50,000$. For each trajectory, we consider the sample lengths $n\in\{5000,10000,\ldots,50000\}$ and compute the plug-in transfer entropy estimator introduced in Section~\ref{sec:plugin-estimators}. Since the joint alphabet is $\mathcal Z=\mathcal X\times\mathcal Y=\{0,1\}^2$, we choose $k(n) = \left\lfloor \frac{\log n}{2\log 4}\right\rfloor.$ Hence $k(n)\to\infty$ and $k(n)\leq C\log n$, with $C=\frac{1}{2\log 4}<\frac{1}{\log 4},$ so that the growth condition of Theorem~\ref{thm:concentration_transfer_entropy_rate} is satisfied.

Figure~\ref{fig:transfer_entropy_convergence} shows the Monte Carlo median of $\widehat T_n^{k(n)}\left(X_{-k(n)+1}^{0}\to Y_{-k(n)+1}^{0}\right),$ together with the empirical $5\%$--$95\%$ quantile band. The dashed horizontal line represents the exact value of $T(\boldsymbol X\to\boldsymbol Y)$. As the sample size increases, the empirical distribution becomes increasingly concentrated around the target value, while the median remains close to the transfer entropy rate. This behavior is consistent with the almost-sure convergence established in Corollary~\ref{cor:concentration_transfer_entropy_rate}.

$\quad$

\begin{figure}[http!]
\centering
\includegraphics[scale=.30]{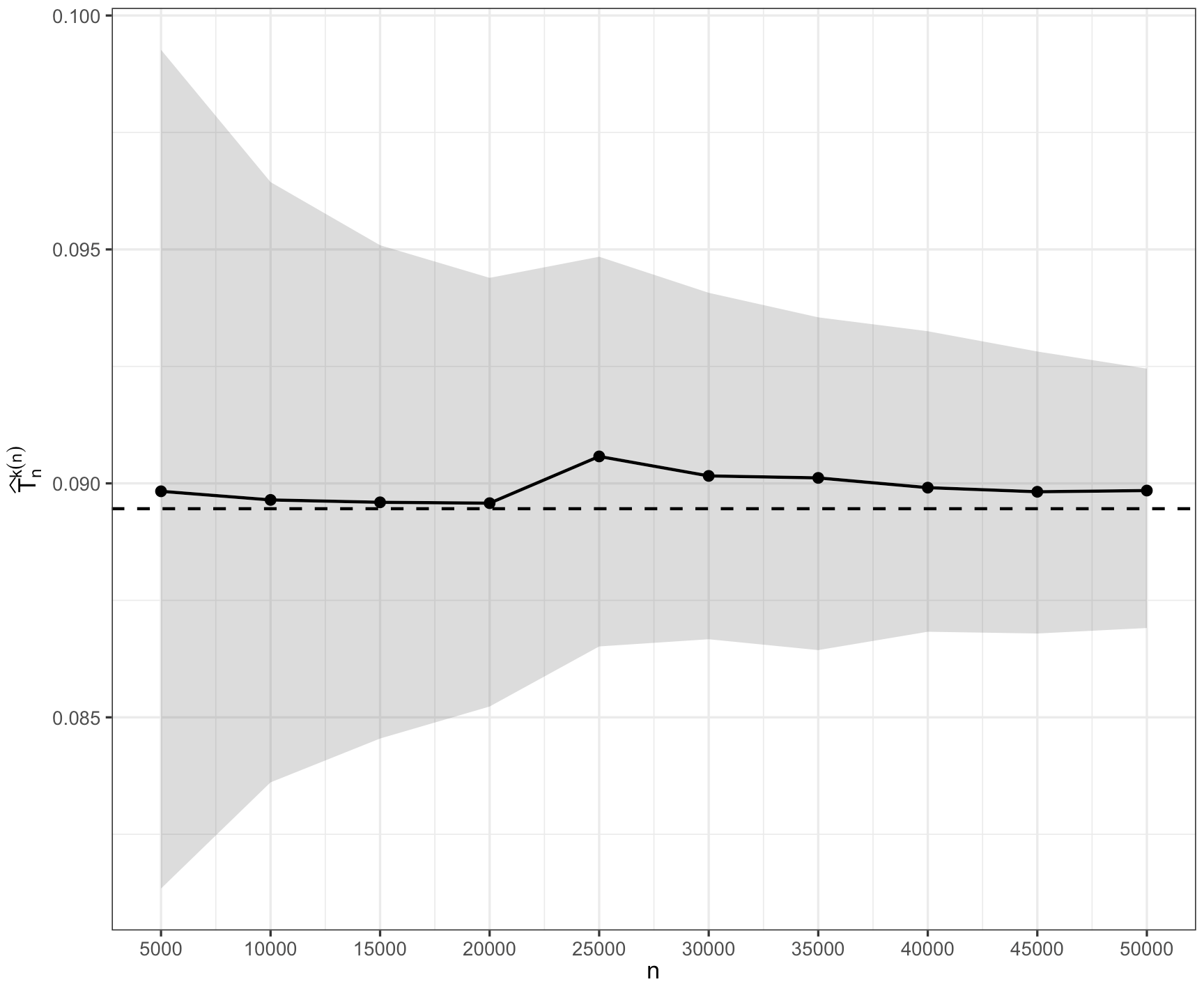}
\caption{Convergence of the plug-in transfer-entropy-rate estimator for the bivariate chain of Example~\ref{ex:transfer-entropy-non-Markovian}.}
\label{fig:transfer_entropy_convergence}
\end{figure}

To examine concentration more directly, Figure~\ref{fig:transfer_entropy_exceedance}
displays the empirical exceedance probabilities
$$
\widehat P_n^T(\varepsilon) := \frac{1}{R}\sum_{r=1}^{R}\mathbf 1\left\{\left|
\widehat T_n^{(r),k(n)} - T(\boldsymbol X\to\boldsymbol Y) \right| > \varepsilon
\right\},
$$
for $\varepsilon\in\{0.002,0.003,0.004\}.$ For each fixed deviation level, the empirical exceedance probabilities show a clear overall decrease as $n$ increases. Moreover, larger values of $\varepsilon$ lead to smaller exceedance probabilities, as expected. This provides a finite-sample illustration of the concentration phenomenon described by Theorem~\ref{thm:concentration_transfer_entropy_rate}.

Finally, Theorem~\ref{thm:concentration_transfer_entropy_rate} gives an exponential upper bound whose relevant scale is {\footnotesize{$e^{-c\,\frac{n}{\log^4 n}.}$}}. Motivated by this rate, Figure~\ref{fig:transfer_entropy_scale} plots $-\log\widehat P_n(\varepsilon)$ against $n/\log^4 n$. The curves exhibit an increasing and approximately linear trend over the range considered, particularly for the larger deviation levels. This behavior is consistent with an exponential-type decay on the scale $n/\log^4 n$ predicted by Theorem~\ref{thm:concentration_transfer_entropy_rate}.

\begin{figure}[http!]
\centering
\includegraphics[scale=.30]{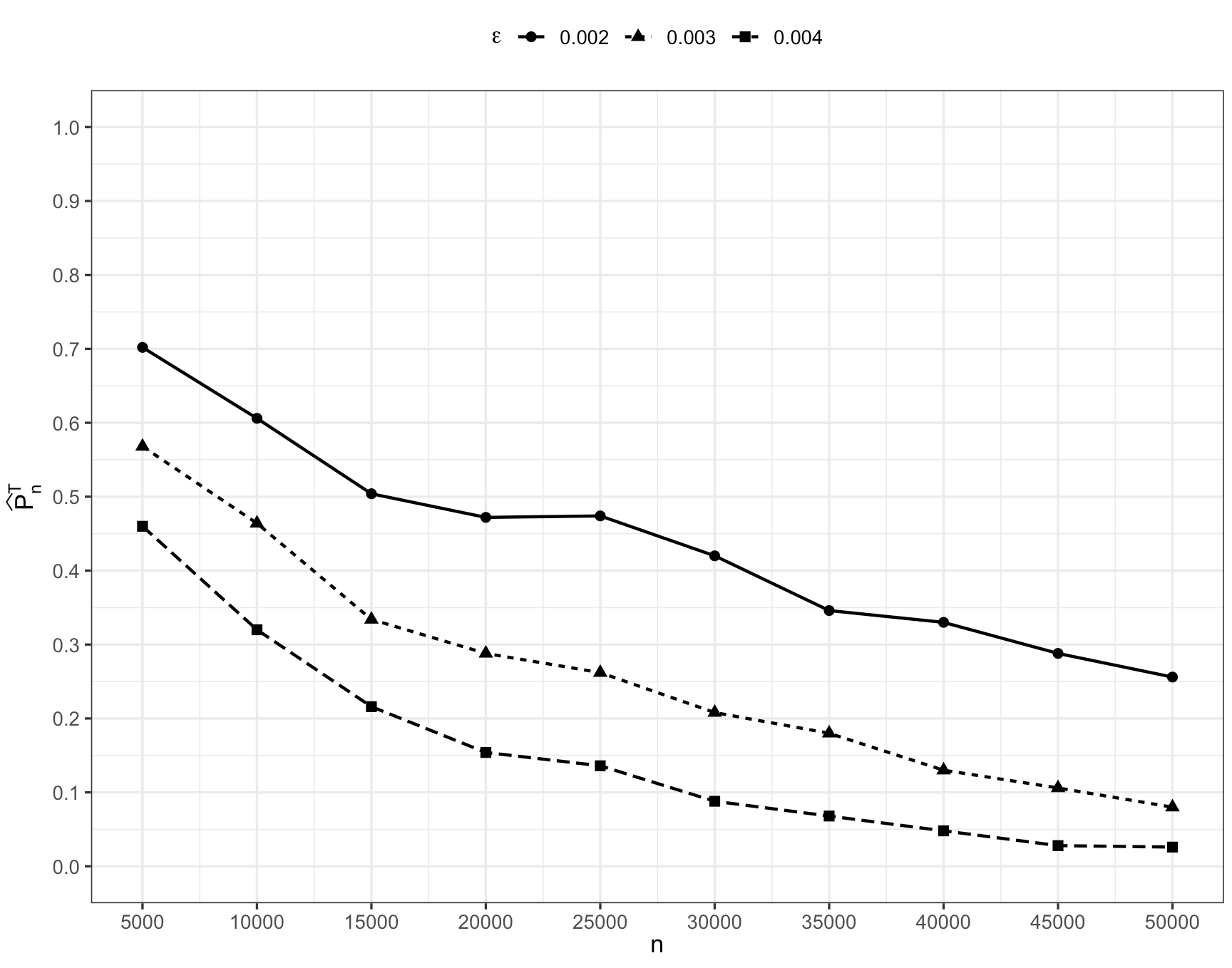}
\caption{Empirical exceedance probabilities of the transfer-entropy-rate estimator for the bivariate chain of Example~\ref{ex:transfer-entropy-non-Markovian}.}
\label{fig:transfer_entropy_exceedance}
\end{figure}

\begin{figure}[http!]
\centering
\includegraphics[scale=.30]{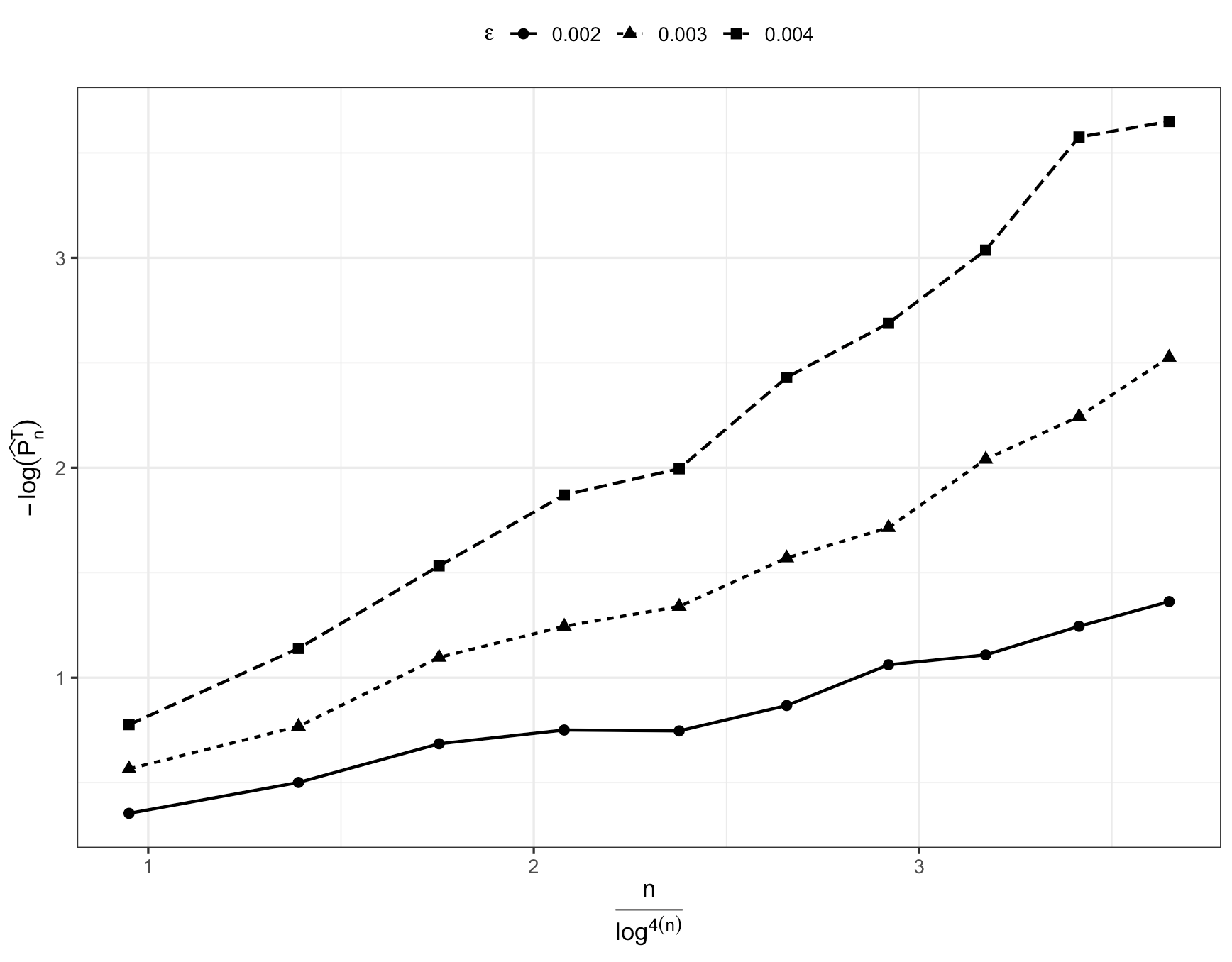}
\caption{Exponential-type concentration scale for the transfer-entropy-rate estimator of the bivariate chain of Example~\ref{ex:transfer-entropy-non-Markovian}.}
\label{fig:transfer_entropy_scale}
\end{figure}

\newpage

\section{Proofs}
\label{sec:proofs}

\subsection{Proof of Theorem \ref{thm:concentration_entropy_mean}}
\label{sec:proof_theorem_concentration_entropy_mean}

Let $y_{-n}^{-1},z_{-n}^{-1}\in\mathcal Y^n$ differ at only one coordinate. Changing one coordinate affects at most $k$ empirical $k$-blocks and at most $k-1$ empirical $(k-1)$-blocks. Since replacing one observation in an empirical distribution based on $n$ observations
changes its entropy by at most $(\log n+1)/n$, we have
\[
\left|\widehat H_{n,k}(y_{-n}^{-1})-\widehat H_{n,k}(z_{-n}^{-1})\right|\leq\frac{k(\log n+1)}{n},
\]
and
\[
\left|\widehat H_{n,k-1}(y_{-n}^{-1})-\widehat H_{n,k-1}(z_{-n}^{-1})\right|\leq\frac{(k-1)(\log n+1)}{n}.
\]

Define
\[
F_{n,k}(Y_{-n}^{-1}) :=\widehat H_n\!\left(Y_0\mid Y_{-k+1}^{-1}\right).
\]
For $i=-n,\ldots,-1$, define the coordinate oscillation of $F_{n,k}$ by
\[
\delta_i(F_{n,k}):=\sup\left\{\left|F_{n,k}(y_{-n}^{-1})-F_{n,k}(z_{-n}^{-1})\right|
:y_j=z_j \text{ for all } j\neq i\right\}.
\]
Therefore, since $F_{n,k} = \widehat H_{n,k} - \widehat H_{n,k-1},$ the coordinate oscillations of $F_{n,k}$ satisfy
\[
\delta_i(F_{n,k}) \leq \frac{(2k-1)(\log n+1)}{n}, \qquad i=-n,\ldots,-1.
\]
Hence,
\begin{equation}
\|\delta F_{n,k}\|_{\ell^2}^2 := \sum_{i=-n}^{-1}\delta_i(F_{n,k})^2 \leq \frac{(2k-1)^2(\log n+1)^2}{n}.
\label{eq:l2-oscillation-F}
\end{equation}

By the concentration inequality for coupling-from-the-past processes
\cite[Corollary 1]{gallo2014attractive}, Assumption~(II) and
\eqref{eq:l2-oscillation-F} imply that, for every $\xi>0$,
\begin{equation}
\mu_Y\left(\left|F_{n,k} - E_{\mu_Y}[F_{n,k}]\right|>\xi\right)
\leq 4\exp\left\{-\frac{2\xi^2 n}{9(1+E[\theta])^2(2k-1)^2(\log n+1)^2}
\right\}.
\label{eq:concentration-around-mean}
\end{equation}
\qed

\subsection{Proof of Remark \ref{rem:entropy_logarithmic_blocks}}
\label{sec:proof_remark_entropy_logarithmic_blocks}

Now take $k=k(n)$ and assume that, for all sufficiently large $n$,
$k(n)\leq C\log n$, for some constant $C>0$. Then, for all sufficiently
large $n$,
\[
(2k(n)-1)^2(\log n+1)^2
\leq
16C^2\log^4 n.
\]
Hence, by \eqref{eq:concentration-around-mean}, for every $\xi>0$ and all
sufficiently large $n$,
\[
\mu_Y\left(
\left|
\widehat H_n^{k(n)}
-
E_{\mu_Y}\!\left[
\widehat H_n^{k(n)}
\right]
\right|
>\xi
\right)
\leq
4\exp\left\{
-\frac{
\xi^2 n
}{
72C^2(1+E[\theta])^2\log^4 n
}
\right\},
\]
where, for shorthand notation, we write $\widehat H_n^{k(n)} := H_n\!\left(Y_0\mid Y_{-k(n)+1}^{-1}\right)$.

\qed

\subsection{Proof of Theorem \ref{thm:concentration_entropy_rate}}
\label{sec:proof_theorem_entropy_rate}

We first derive an upper bound for the bias
$\left| E_{\mu_Y}\!\left[\widehat H_n\!\left(Y_0\mid Y_{-k+1}^{-1}\right)\right]-h(\boldsymbol Y)\right|.$
Fix integers $k$ and $n$ such that $2\leq k\leq n$. For
$b_{-k+1}^0\in\mathcal Y^k$, write
$\mu_k(b_{-k+1}^0) := \mu_Y([b_{-k+1}^0])$, and recall the shorthand
notation $\widehat\mu_{n,k}(b_{-k+1}^0)$ for the empirical $k$-block
distribution introduced in Section~\ref{sec:plugin-estimators}.

We also define the \emph{finite-past transition kernel}
$g_k: \mathcal{Y}^{k-1} \times \mathcal{Y} \rightarrow [0,1]$ such that
\[
g_k\!\left(b_0\mid b_{-k+1}^{-1}\right):=
\frac{\mu_Y([b_{-k+1}^0])}{\mu_Y([b_{-k+1}^{-1}])}.
\]
By Assumption~(I) and compatibility of $\mu_Y$ with $g$,
\[
\begin{aligned}
\mu_Y([b_{-k+1}^0])
&=
\int_{\mathcal Y^{\mathbb Z}}
\mathbf 1_{[b_{-k+1}^{-1}]}(y)
\mathbf 1_{[b_0]}(y)\,d\mu_Y(y)\\
&=
\int_{[b_{-k+1}^{-1}]}
E_{\mu_Y}\!\left[
\mathbf 1_{[b_0]}\,\middle|\,\mathcal F^{\leq -1}
\right](y)\,d\mu_Y(y)\\
&=
\int_{[b_{-k+1}^{-1}]}
g\!\left(b_0\mid y_{-\infty}^{-1}\right)\,d\mu_Y(y)\\
&\geq
\epsilon_0\,\mu_Y([b_{-k+1}^{-1}]).
\end{aligned}
\]
Iterating this inequality shows that every finite cylinder has a positive probability. In particular,
\[
g_k\!\left(b_0\mid b_{-k+1}^{-1}\right)
=
\frac{\mu_Y([b_{-k+1}^0])}{\mu_Y([b_{-k+1}^{-1}])}
\geq \epsilon_0.
\]

Let
\[
D(\nu\Vert\rho) := \sum_i \nu(i)\log\frac{\nu(i)}{\rho(i)}
\]
denote the \emph{relative entropy} of $\nu$ with respect to $\rho$. An
important consequence of the periodic construction is that the empirical
block distributions are consistent under marginalization. Indeed, for every
$b_{-k+1}^{-1}\in\mathcal Y^{k-1}$,
\[
\sum_{b_0\in\mathcal Y}
\widehat\mu_{n,k}(b_{-k+1}^0)
=
\widehat\mu_{n,k-1}(b_{-k+1}^{-1}).
\]
Therefore,
\begin{align}
\widehat H_n\!\left(Y_0\mid Y_{-k+1}^{-1}\right)
&=
-\sum_{b_{-k+1}^0\in\mathcal Y^k}
\widehat\mu_{n,k}(b_{-k+1}^0)
\log
\frac{\widehat\mu_{n,k}(b_{-k+1}^0)}
{\widehat\mu_{n,k-1}(b_{-k+1}^{-1})}
\nonumber\\
&=
-\sum_{b_{-k+1}^0\in\mathcal Y^k}
\widehat\mu_{n,k}(b_{-k+1}^0)
\log g_k\!\left(b_0\mid b_{-k+1}^{-1}\right)
+
\Delta\widehat D_{n,k},
\label{eq:bias-KL-decomposition}
\end{align}
where
\[
\Delta\widehat D_{n,k}
:=
-D\!\left(\widehat\mu_{n,k}\Vert\mu_k\right)
+
D\!\left(\widehat\mu_{n,k-1}\Vert\mu_{k-1}\right).
\]

We now take expectation in \eqref{eq:bias-KL-decomposition}. Because
$\widehat\mu_{n,k}$ is defined from the periodic extension of the observed
sample, its expectation is not, in general, exactly equal to $\mu_k$. To
isolate this boundary effect, define
\[
R_{n,k}
:=
-\sum_{b_{-k+1}^0\in\mathcal Y^k}
\left(
E_{\mu_Y}\!\left[\widehat\mu_{n,k}(b_{-k+1}^0)\right]
-
\mu_k(b_{-k+1}^0)
\right)
\log g_k\!\left(b_0\mid b_{-k+1}^{-1}\right).
\]
It follows that
\begin{align}
E_{\mu_Y}\!\left[
\widehat H_n\!\left(Y_0\mid Y_{-k+1}^{-1}\right)
\right]
&=
H\!\left(Y_0\mid Y_{-k+1}^{-1}\right)
+
R_{n,k}
+
E_{\mu_Y}\!\left[\Delta\widehat D_{n,k}\right].
\label{eq:expected-empirical-entropy}
\end{align}

We next control the term $R_{n,k}$. Recall from
Section~\ref{sec:plugin-estimators} the auxiliary empirical distribution
defined in \eqref{eq:auxiliary_empirical_distribution}, i.e.,
\[
\overline\mu_{n,k}(b_{-k+1}^0)
:=
\frac{1}{n}
\sum_{i=-n}^{-1}
\mathbf 1\left\{
Y_i^{i+k-1}=b_{-k+1}^0
\right\}.
\]
By stationarity,
\[
E_{\mu_Y}\!\left[
\overline\mu_{n,k}(b_{-k+1}^0)
\right]
=
\mu_k(b_{-k+1}^0).
\]
Moreover, $\widehat\mu_{n,k}$ and $\overline\mu_{n,k}$ may differ only
at the last $k-1$ starting positions. Hence, for every realization,
\[
\sum_{b_{-k+1}^0\in\mathcal Y^k}
\left|
\widehat\mu_{n,k}(b_{-k+1}^0)
-
\overline\mu_{n,k}(b_{-k+1}^0)
\right|
\leq
\frac{2(k-1)}{n}.
\]
Consequently,
\begin{equation}
\left\|
E_{\mu_Y}[\widehat\mu_{n,k}]
-
\mu_k
\right\|_{\mathrm{TV}}
=
\dfrac{1}{2}
\sum_{b_{-k+1}^0\in\mathcal Y^k}
\left|
E_{\mu_Y}\!\left[\widehat\mu_{n,k}(b_{-k+1}^0)\right]
-
\mu_k(b_{-k+1}^0)
\right|
\leq
\frac{k-1}{n}.
\label{eq:periodic-boundary-l1}
\end{equation}
By Assumption~(I),
\[
0
\leq
-\log g_k\!\left(b_0\mid b_{-k+1}^{-1}\right)
\leq
\log\frac{1}{\epsilon_0}.
\]
Since both $E_{\mu_Y}[\widehat\mu_{n,k}]$ and $\mu_k$ are probability
distributions, their difference has total mass zero. Combining these facts
with \eqref{eq:periodic-boundary-l1} gives
\begin{equation}
|R_{n,k}|
\leq
\frac{k-1}{n}
\log\frac{1}{\epsilon_0}.
\label{eq:periodic-boundary-error}
\end{equation}

It remains to identify the error produced by replacing the infinite past by
a finite past. We define
\begin{equation}
\gamma_k(g)
:=
H\!\left(Y_0\mid Y_{-k+1}^{-1}\right)
-
h(\boldsymbol Y).
\label{eq:def-gamma-k}
\end{equation}
By Section~\ref{sec:entropy_quantities}, we know that
$\gamma_k(g)\geq0$ and $\gamma_k(g)\longrightarrow0$ as $k\to\infty$.
So, subtracting $h(\boldsymbol Y)$ from
\eqref{eq:expected-empirical-entropy} and using
\eqref{eq:def-gamma-k}, we obtain
\[
E_{\mu_Y}\!\left[
\widehat H_n\!\left(Y_0\mid Y_{-k+1}^{-1}\right)
\right]
-
h(\boldsymbol Y)
=
\gamma_k(g)
+
R_{n,k}
+
E_{\mu_Y}\!\left[\Delta\widehat D_{n,k}\right].
\]
Therefore, by \eqref{eq:periodic-boundary-error},
\begin{equation}
\left|
E_{\mu_Y}\!\left[
\widehat H_n\!\left(Y_0\mid Y_{-k+1}^{-1}\right)
\right]
-
h(\boldsymbol Y)
\right|
\leq
\gamma_k(g)
+
\frac{k-1}{n}\log\frac{1}{\epsilon_0}
+
\left|
E_{\mu_Y}\!\left[\Delta\widehat D_{n,k}\right]
\right|.
\label{eq:bias-decomposition}
\end{equation}

We next control the relative-entropy term appearing in
\eqref{eq:bias-decomposition}. Recall that
\[
\Delta\widehat D_{n,k}
=
-D\!\left(\widehat\mu_{n,k}\Vert\mu_k\right)
+
D\!\left(\widehat\mu_{n,k-1}\Vert\mu_{k-1}\right).
\]
Since $\widehat\mu_{n,k-1}$ and $\mu_{k-1}$ are, respectively, the
$(k-1)$-block marginals of $\widehat\mu_{n,k}$ and $\mu_k$, the
data-processing inequality for relative entropy \citep{Csiszar1967} yields
\[
D\!\left(
\widehat\mu_{n,k-1}\Vert\mu_{k-1}
\right)
\leq
D\!\left(
\widehat\mu_{n,k}\Vert\mu_k
\right),
\]
which implies that $\Delta\widehat D_{n,k}\leq0$, and therefore
\begin{equation}
\left|
E_{\mu_Y}\!\left[
\Delta\widehat D_{n,k}
\right]
\right|
\leq
E_{\mu_Y}\!\left[
D\!\left(\widehat\mu_{n,k}\Vert\mu_k\right)
\right].
\label{eq:delta-D-control}
\end{equation}

The empirical distribution $\widehat\mu_{n,k}$ is defined through the
periodic extension of the sample. Since the mixing assumption applies to the
original process rather than to this periodic extension, we separate the
blocks entirely contained in the observed sample from those affected by
periodicization. For $b_{-k+1}^{0}\in\mathcal Y^k$, define
\[
\widehat\nu_{n,k}(b_{-k+1}^{0})
:=
\frac{1}{n-k+1}
\sum_{i=-n}^{-k}
\mathbf 1\left\{
Y_i^{i+k-1}=b_{-k+1}^{0}
\right\}.
\]

Thus, $\hat{\nu}_{n,k}$ is the empirical distribution of the $n-k+1$
blocks with length $k$ entirely contained in $Y_{-n}^{-1}$. The remaining
$k-1$ blocks are those completed using the periodic extension, whose
empirical distribution $\hat{\rho}_{n,k}$ is such that
\begin{equation}
\widehat\mu_{n,k}
=
\frac{n-k+1}{n}\hat{\nu}_{n,k}
+
\frac{k-1}{n}\hat{\rho}_{n,k}.
\label{eq:periodic-mixture}
\end{equation}

By convexity of relative entropy,
\begin{equation}
D\!\left(\widehat\mu_{n,k}\Vert\mu_k\right)
\leq
\frac{n-k+1}{n}
D\!\left(\hat{\nu}_{n,k}\Vert\mu_k\right)
+
\frac{k-1}{n}
D\!\left(\hat{\rho}_{n,k}\Vert\mu_k\right).
\label{eq:KL-convexity}
\end{equation}
Consequently, for any probability distribution $\rho$ on $\mathcal Y^k$,
\[
D(\rho\Vert\mu_k)
=
\sum_{b_{-k+1}^{0}\in\mathcal Y^k}
\rho(b_{-k+1}^{0})
\log
\frac{\rho(b_{-k+1}^{0})}
{\mu_k(b_{-k+1}^{0})}
\leq
-\sum_{b_{-k+1}^{0}\in\mathcal Y^k}
\rho(b_{-k+1}^{0})
\log\mu_k(b_{-k+1}^{0})
\leq
k\log\frac{1}{\epsilon_0},
\]
where the last inequality follows from Assumption~(I). Applying this
inequality to $\hat{\rho}_{n,k}$ in \eqref{eq:KL-convexity} gives
\begin{equation}
E_{\mu_Y}\!\left[
D\!\left(\widehat\mu_{n,k}\Vert\mu_k\right)
\right]
\leq
\frac{n-k+1}{n}
E_{\mu_Y}\!\left[
D\!\left(\hat{\nu}_{n,k}\Vert\mu_k\right)
\right]
+
\frac{k(k-1)}{n}
\log\frac{1}{\epsilon_0}.
\label{eq:periodic-KL-control}
\end{equation}

Relative entropy is bounded above by the chi-square divergence
\cite[Lemma~2.7]{Tsybakov2009}, so that
\[
D(\hat{\nu}_{n,k}\Vert\mu_k)
\leq
\sum_{b_{-k+1}^{0}\in\mathcal Y^k}
\frac{
\left(
\hat{\nu}_{n,k}(b_{-k+1}^{0})
-
\mu_k(b_{-k+1}^{0})
\right)^2
}{
\mu_k(b_{-k+1}^{0})
}.
\]
Taking expectations and using stationarity, we have
\begin{equation}
E_{\mu_Y}\!\left[
D(\hat{\nu}_{n,k}\Vert\mu_k)
\right]
\leq
\sum_{b_{-k+1}^{0}\in\mathcal Y^k}
\frac{
\operatorname{Var}_{\mu_Y}
\left(
\hat{\nu}_{n,k}(b_{-k+1}^{0})
\right)
}{
\mu_k(b_{-k+1}^{0})
}.
\label{eq:KL-variance}
\end{equation}

To bound these variances, fix $b_{-k+1}^{0}\in\mathcal Y^k$ and set
$I_i^{b_{-k+1}^{0}}
:=
\mathbf 1\left\{
Y_i^{i+k-1}=b_{-k+1}^{0}
\right\}.$
By stationarity,
\[
\begin{aligned}
\operatorname{Var}_{\mu_Y}
\left(
\hat{\nu}_{n,k}(b_{-k+1}^{0})
\right)
\leq{}&
\frac{\mu_k(b_{-k+1}^{0})}{n-k+1}\\
&+
\frac{2}{(n-k+1)^2}
\sum_{r=1}^{n-k}
(n-k+1-r)
\left|
\operatorname{Cov}_{\mu_Y}
\left(
I_0^{b_{-k+1}^{0}},
I_r^{b_{-k+1}^{0}}
\right)
\right|.
\end{aligned}
\]
For $1\leq r<k$, the two blocks overlap, and
\[
\left|
\operatorname{Cov}_{\mu_Y}
\left(
I_0^{b_{-k+1}^{0}},
I_r^{b_{-k+1}^{0}}
\right)
\right|
\leq
\mu_k(b_{-k+1}^{0}).
\]
If $r\geq k$, then
$I_0^{b_{-k+1}^{0}}$ and $I_r^{b_{-k+1}^{0}}$ are measurable with respect
to sigma-fields separated by a lag $r-k+1$. Hence, by Assumption~(III),
\[
\left|
\operatorname{Cov}_{\mu_Y}
\left(
I_0^{b_{-k+1}^{0}},
I_r^{b_{-k+1}^{0}}
\right)
\right|
\leq
\beta(r-k+1)
\leq
\gamma e^{-\delta(r-k+1)}.
\]
Therefore,
\[
\left|
\operatorname{Cov}_{\mu_Y}
\left(
I_0^{b_{-k+1}^{0}},
I_r^{b_{-k+1}^{0}}
\right)
\right|
\leq
\min\left\{
\mu_k(b_{-k+1}^{0}),
\gamma e^{-\delta(r-k+1)}
\right\},
\qquad r\geq k.
\]
It follows that
\begin{equation}
\begin{aligned}
\operatorname{Var}_{\mu_Y}
\left(
\hat{\nu}_{n,k}(b_{-k+1}^{0})
\right)
\leq
\frac{\mu_k(b_{-k+1}^{0})}{n-k+1}
\Bigg[
2k-1
+
2\sum_{s=1}^{\infty}
\min
\left\{
1,
\frac{\gamma e^{-\delta s}}
{\mu_k(b_{-k+1}^{0})}
\right\}
\Bigg].
\end{aligned}
\label{eq:variance-nu}
\end{equation}

Define
\[
A_k
:=
\sum_{s=1}^{\infty}
\min\left\{
1,
\gamma\epsilon_0^{-k}e^{-\delta s}
\right\}.
\]
By \eqref{eq:variance-nu} and since
$\mu_k(b_{-k+1}^{0})\geq\epsilon_0^k$, we have
\begin{equation}
\frac{
\operatorname{Var}_{\mu_Y}
\left(
\hat{\nu}_{n,k}(b_{-k+1}^{0})
\right)
}{
\mu_k(b_{-k+1}^{0})
}
\leq
\frac{2k-1+2A_k}{n-k+1}.
\label{eq:variance-ratio}
\end{equation}
Splitting the geometric series at the index where
$\gamma\epsilon_0^{-k}e^{-\delta s}$ becomes smaller than one gives
\begin{equation}
A_k
\leq
1+
\frac{\log_+(\gamma\epsilon_0^{-k})}{\delta}
+
\frac{1}{e^\delta-1},
\label{eq:Ak-bound}
\end{equation}
where $\log_+(x):=\max\{\log x,0\}$. In particular,
\[
A_k
\leq
1+
\frac{\log_+\gamma}{\delta}
+
\frac{k}{\delta}\log\frac{1}{\epsilon_0}
+
\frac{1}{e^\delta-1}.
\]
Thus, setting
\[
c_1
:=
2+\frac{2}{\delta}\log\frac{1}{\epsilon_0}
\qquad\text{and}\qquad
c_0
:=
1+
\frac{2\log_+\gamma}{\delta}
+
\frac{2}{e^\delta-1},
\]
we have
\[
2k-1+2A_k
\leq
c_1k+c_0.
\]
Combining this inequality with \eqref{eq:KL-variance} and
\eqref{eq:variance-ratio}, and summing over
$b_{-k+1}^{0}\in\mathcal Y^k$, we obtain
\begin{equation}
E_{\mu_Y}\!\left[
D(\hat{\nu}_{n,k}\Vert\mu_k)
\right]
\leq
\frac{|\mathcal Y|^k}{n-k+1}
(c_1k+c_0).
\label{eq:KL-nu-final}
\end{equation}

Substituting \eqref{eq:KL-nu-final} into
\eqref{eq:periodic-KL-control}, and then using
\eqref{eq:delta-D-control}, yields
\begin{equation}
\left|
E_{\mu_Y}\!\left[
\Delta\widehat D_{n,k}
\right]
\right|
\leq
\frac{|\mathcal Y|^k}{n}
(c_1k+c_0)
+
\frac{k(k-1)}{n}
\log\frac{1}{\epsilon_0}.
\label{eq:Delta-D-final}
\end{equation}

Finally, inserting \eqref{eq:Delta-D-final} into
\eqref{eq:bias-decomposition}, we conclude that
\begin{equation}
\left|
E_{\mu_Y}\!\left[
\widehat H_n\!\left(Y_0\mid Y_{-k+1}^{-1}\right)
\right]
-
h(\boldsymbol Y)
\right|
\leq
\gamma_k(g)
+
\frac{|\mathcal Y|^k}{n}(c_1k+c_0)
+
\frac{k^2-1}{n}\log\frac{1}{\epsilon_0}.
\label{eq:final-bias-bound}
\end{equation}

We now take $k=k(n)$ and define
\[
r_n
:=
\gamma_{k(n)}(g)
+
\frac{|\mathcal Y|^{k(n)}}{n}
\bigl(c_1k(n)+c_0\bigr)
+
\frac{k(n)^2-1}{n}
\log\frac{1}{\epsilon_0}.
\]
By \eqref{eq:final-bias-bound},
\[
\left|
E_{\mu_Y}\!\left[
\widehat H_n\!\left(
Y_0\mid Y_{-k(n)+1}^{-1}
\right)
\right]
-
h(\boldsymbol Y)
\right|
\leq
r_n.
\]
Since $k(n)\to\infty$, $k(n)\leq C\log n$, and
$0<C<1/\log|\mathcal Y|$, we have $r_n\to0$.

By the triangle inequality,
\[
\left\{
\left|
\widehat H_n^{k(n)}
-
h(\boldsymbol Y)
\right|
>
\xi+r_n
\right\}
\subseteq
\left\{
\left|
\widehat H_n^{k(n)}
-
E_{\mu_Y}\!\left[
\widehat H_n^{k(n)}
\right]
\right|
>
\xi
\right\}.
\]
Therefore, by Remark~\ref{rem:entropy_logarithmic_blocks}, for every
$\xi>0$ and all sufficiently large $n$,
\[
\mu_Y\Bigg(
\left|
\widehat H_n\!\left(
Y_0\mid Y_{-k(n)+1}^{-1}
\right)
-
h(\boldsymbol Y)
\right|
>
\xi+r_n
\Bigg)
\leq
4\exp\left\{
-\frac{
\xi^2 n
}{
72C^2(1+E[\theta])^2\log^4 n
}
\right\}.
\]
\qed

\subsection{Proof of Corollary \ref{cor:concentration_entropy_rate}}
\label{sec:proof_corollary_entropy_rate}

Since $r_n\to0$ and
\[
\sum_{n\ge2}\exp\left\{-c\,\frac{n}{\log^4 n}\right\}<\infty,
\]
with $c > 0$, Theorem~\ref{thm:concentration_entropy_rate} and the Borel--Cantelli lemma \cite[Theorem~2.3.1]{Durrett2019} imply
\[
\widehat H_n\!\left(
Y_0\mid Y_{-k(n)+1}^{-1}
\right)
\longrightarrow
h(\boldsymbol Y)
\qquad \mu_Y\text{-a.s.}
\]
\qed

\subsection{Proof of Theorem \ref{thm:concentration_transfer_entropy_mean}}
\label{sec:proof_theorem_te_mean}

Fix integers $k$ and $n$ such that $2\leq k\leq n$. By the definition of the plug-in transfer entropy,
{\footnotesize{
\[
\begin{aligned}\widehat T_n^k\left(X_{-k+1}^{0}\to Y_{-k+1}^{0}\right) &=
\widehat H_n\left(Y_0\mid Y_{-k+1}^{-1}\right) - \widehat H_n\left(Y_0\mid X_{-k+1}^{-1},Y_{-k+1}^{-1}\right)\\
&=
\widehat H_n\left(Y_{-k+1}^{0}\right)-\widehat H_n\left(Y_{-k+1}^{-1}\right)-\widehat H_n\left(X_{-k+1}^{-1},Y_{-k+1}^{0}\right)+\widehat H_n\left(X_{-k+1}^{-1},Y_{-k+1}^{-1}\right).
\end{aligned}
\]}}

Let $(x_{-n}^{-1},y_{-n}^{-1})$ and $(x_{-n}^{\prime\,-1},y_{-n}^{\prime\,-1})$ be two samples that differ only at coordinate $i\in\{-n,\ldots,-1\}$. Changing the pair $(x_i,y_i)$ affects at most $k$ empirical blocks associated with $\widehat H_n\left(Y_{-k+1}^{0}\right)$ and 
$\widehat H_n\left(X_{-k+1}^{-1},Y_{-k+1}^{0}\right)$, and at most $k-1$ empirical blocks associated with $\widehat H_n\left(Y_{-k+1}^{-1}\right)$ and  $\widehat H_n\left(X_{-k+1}^{-1},Y_{-k+1}^{-1}\right).$  Since replacing one observation in an empirical distribution based on $n$ observations changes its entropy by at most $(\log n+1)/n$, applying this
bound successively to the affected empirical blocks shows that the absolute differences between the corresponding empirical entropies computed from the two samples are bounded by $\frac{k(\log n+1)}{n}$ and $\frac{(k-1)(\log n+1)}{n},$ respectively.

For shorthand notation, write $\widehat T_n^k := \widehat T_n^k
\left(X_{-k+1}^{0}\to Y_{-k+1}^{0}\right)$. For each $i\in\{-n,\ldots,-1\}$, define the coordinate oscillation of the plug-in transfer entropy estimator by
\[
\delta_i\!\left(
\widehat T_n^k
\right) := \sup\Bigg\{\left| \widehat T_n^k - \widehat T_n^{\prime\,k} \right| :
(x_j,y_j)=(x_j',y_j') \ \text{for all } j\neq i \Bigg\},
\]
where  $\widehat T_n^k$ and $\widehat T_n^{\prime\,k}$ are empirical transfer entropies computed from $(x_{-n}^{-1},y_{-n}^{-1})$ and $(x_{-n}^{\prime\,-1},y_{-n}^{\prime\,-1})$, respectively. Therefore, for each $i\in\{-n,\ldots,-1\}$, the triangle inequality gives
\[
\delta_i\!\left(\widehat T_n^k\right)
\leq
\frac{2k(\log n+1)}{n}
+
\frac{2(k-1)(\log n+1)}{n},
\]
and hence
\[
\delta_i\!\left(\widehat T_n^k\right)
\leq
\frac{(4k-2)(\log n+1)}{n}.
\]
Consequently,
\begin{equation}
\left\|\delta \widehat T_n^k\right\|_{\ell^2}^2 := \sum_{i=-n}^{-1} \delta_i \!\left(\widehat T_n^k\right)^2\leq\frac{(4k-2)^2(\log n+1)^2}{n}.
\label{eq:l2-oscillation-T}
\end{equation}

By the concentration inequality for coupling-from-the-past processes \cite[Corollary 1]{gallo2014attractive}, Assumption~(II) and \eqref{eq:l2-oscillation-T} imply that, for every $\xi>0$,
\begin{equation}
\mu\left(\left|\widehat T_n^k-E_\mu\!\left[\widehat T_n^k\right]\right|>\xi\right)\leq
4\exp\left\{-\frac{2\xi^2 n}{9(1+E[\theta])^2(4k-2)^2(\log n+1)^2}\right\}.
\label{eq:concentration_TE_1}
\end{equation}

\qed

\subsection{Proof of Remark \ref{rmk:concentration_transfer_entropy_mean}}
\label{sec:proof_remark_te_mean}

Now take $k=k(n)$ and assume that, for all sufficiently large $n$, $k(n)\leq C\log n,$ for some constant $C>0$. Then, for all sufficiently large $n$, 
$$(4k(n)-2)^2(\log n+1)^2 \leq 64C^2\log^4 n.$$
Hence, by \eqref{eq:concentration_TE_1}, for every $\xi>0$ and all sufficiently large $n$, 
\[
\mu\left(\left|\widehat T_n^{k(n)}-E_\mu\!\left[\widehat T_n^{k(n)}\right]\right|>\xi\right)
\leq 4\exp\left\{-\frac{\xi^2 n}{288C^2(1+E[\theta])^2\log^4 n}\right\}.
\]

\qed

\subsection{Proof of Theorem \ref{thm:concentration_transfer_entropy_rate}}
\label{sec:proof_theorem_te_rate}

Fix integers $k$ and $n$ such that $2\leq k\leq n$. We first derive an upper bound for the bias
$\left|E_\mu\!\left[\widehat T_n^k\left(X_{-k+1}^{0}\to Y_{-k+1}^{0}\right)\right] -
T(\boldsymbol X\to\boldsymbol Y)\right|.$ By shorthand notation, write $\widehat T_n^k := \widehat T_n^k\left(X_{-k+1}^{0}\to Y_{-k+1}^{0}\right)$. By the definitions of the plug-in transfer entropy and the transfer entropy rate,
\begin{equation}
\left|
E_\mu\!\left[\widehat T_n^k\right] - T(\boldsymbol X\to\boldsymbol Y) \right| \leq
\gamma_k^{Y}+\gamma_k^{Y|X} + B_{n,k}^{Y} + B_{n,k}^{Y|X}.
\label{eq:TE-bias-decomposition}
\end{equation}
where 
\begin{align*}
\gamma_k^{Y} &:= H\left(Y_0\mid Y_{-k+1}^{-1}\right)-h(\boldsymbol Y), \\
\gamma_k^{Y|X} &:= H\left(Y_0\mid X_{-k+1}^{-1},Y_{-k+1}^{-1}\right) - h(\boldsymbol Y\mid\boldsymbol X), \\
B_{n,k}^{Y} &:= \left|E_\mu\!\left[\widehat H_n\!\left(Y_0\mid Y_{-k+1}^{-1}\right)\right] - H\!\left( Y_0\mid Y_{-k+1}^{-1}\right)\right|,\\
B_{n,k}^{Y|X} &:= \left|E_\mu\!\left[\widehat H_n\!\left(Y_0\mid X_{-k+1}^{-1},Y_{-k+1}^{-1}
\right)\right]-H\!\left(Y_0\mid X_{-k+1}^{-1},Y_{-k+1}^{-1}\right)\right|.
\end{align*}

We first control $B_{n,k}^{Y}$. Since $\boldsymbol Y$ is a marginal of $\boldsymbol Z$, the $\beta$-mixing coefficients of $\boldsymbol Y$ and $\boldsymbol Z$, denoted by $\beta_Y$ and $\beta_Z$, respectively, satisfy, for every $r \geq 1$,
\[
\beta_Y(r)\leq \beta_Z(r)\leq \gamma e^{-\delta r},
\]
for some positive constants $\delta$ and $\gamma$. Moreover, by Assumption~{\rm (I)}, for every $b_0\in\mathcal Y$ and every $(a_{-\infty}^{-1},b_{-\infty}^{-1}) \in\mathcal X^{\leq-1}\times\mathcal Y^{\leq-1}$,
\[
g_{Y|X}\left(b_0\mid a_{-\infty}^{-1},b_{-\infty}^{-1}\right) = \sum_{a_0\in\mathcal X}g\left(
(a_0,b_0)\mid(a_{-\infty}^{-1},b_{-\infty}^{-1})\right) \geq |\mathcal X|\epsilon_0 \geq
\epsilon_0.
\]
Hence, by the tower property, $g_Y\left(b_0\mid b_{-\infty}^{-1}\right)\geq\epsilon_0$, $\mu_Y\text{-a.s.}$ Therefore, the same argument used in the proof of Theorem~\ref{thm:concentration_entropy_rate}, in particular the bounds leading to \eqref{eq:periodic-boundary-error} and \eqref{eq:Delta-D-final}, gives
\begin{equation}
B_{n,k}^{Y}
\leq
\frac{|\mathcal Y|^k}{n}(c_1k+c_0)
+
\frac{k^2-1}{n}\log\frac1{\epsilon_0},
\label{eq:BY-bound}
\end{equation}
where $c_1$ and $c_0$ are the constants introduced in the proof of Theorem \ref{thm:concentration_entropy_rate} (see Section~\ref{sec:proof_theorem_entropy_rate}).

We next control $B_{n,k}^{Y|X}$. Let $\mu_{k,k}$ denote the stationary joint $k$-block distribution of $\boldsymbol Z$, and let $\mu_{k-1,k}$ and $\mu_{k-1,k-1}$ denote its corresponding marginals, in analogy with the empirical distributions introduced in Section~2.5.
We define the finite-past conditional kernel by
\[
g_k^{Y|X}\left(b_0\mid a_{-k+1}^{-1},b_{-k+1}^{-1}\right) := \frac{\mu_{k-1,k}\left(
a_{-k+1}^{-1},b_{-k+1}^{0}\right)}{\mu_{k-1,k-1}\left(a_{-k+1}^{-1},b_{-k+1}^{-1}\right)
}.
\]
By Assumption~{\rm (I)} and the tower property, $g_k^{Y|X}\left(b_0\mid a_{-k+1}^{-1},b_{-k+1}^{-1}\right)\geq\epsilon_0.$ Moreover, by the consistency of the periodic empirical distributions under marginalization,
{\footnotesize
\begin{equation}
\widehat H_n^{Y|X}
=
-\sum_{\substack{a_{-k+1}^{-1}\in\mathcal X^{k-1}\\b_{-k+1}^{0}\in\mathcal Y^k}}
\widehat\mu_{n;k-1,k}\left(a_{-k+1}^{-1},b_{-k+1}^{0}\right)
\log g_k^{Y|X}\left(b_0\mid a_{-k+1}^{-1},b_{-k+1}^{-1}\right)
+
\Delta\widehat D_{n,k}^{Y|X},
\label{eq:rewrite_empirical_conditional_entropy}
\end{equation}
}
where $\widehat H_n^{Y|X} := \widehat H_n\left(Y_0\mid X_{-k+1}^{-1},Y_{-k+1}^{-1}\right)$ and
\[
\Delta\widehat D_{n,k}^{Y|X}
:=
-D\left(\widehat\mu_{n;k-1,k}\Vert\mu_{k-1,k}\right)
+
D\left(\widehat\mu_{n;k-1,k-1}\Vert\mu_{k-1,k-1}\right).
\]

Taking expectation in \eqref{eq:rewrite_empirical_conditional_entropy}, we have
\begin{equation}
E_\mu\!\left[\widehat H_n^{Y|X}\right]
=
H\left(Y_0\mid X_{-k+1}^{-1},Y_{-k+1}^{-1}\right)
+
R_{n,k}^{Y|X}
+
E_\mu\!\left[\Delta\widehat D_{n,k}^{Y|X}\right],
\label{eq:EYgivenX-decomposition}
\end{equation}
where
{\footnotesize
\[
R_{n,k}^{Y|X}
:=
-\sum_{\substack{a_{-k+1}^{-1}\in\mathcal X^{k-1}\\b_{-k+1}^{0}\in\mathcal Y^k}}
\left(
E_\mu\!\left[\widehat\mu_{n;k-1,k}\left(a_{-k+1}^{-1},b_{-k+1}^{0}\right)\right]
-
\mu_{k-1,k}\left(a_{-k+1}^{-1},b_{-k+1}^{0}\right)
\right)
\log g_k^{Y|X}\left(b_0\mid a_{-k+1}^{-1},b_{-k+1}^{-1}\right).
\]
}
Then, by \eqref{eq:EYgivenX-decomposition},
\begin{equation}
B_{n,k}^{Y|X}
\leq
\left|R_{n,k}^{Y|X}\right|
+
\left|
E_\mu\!\left[\Delta\widehat D_{n,k}^{Y|X}\right]
\right|.
\label{eq:BYX-boundary-KL}
\end{equation}

To control $R_{n,k}^{Y|X}$, define
\[
\overline\mu_{n;k-1,k}\left(a_{-k+1}^{-1},b_{-k+1}^{0}\right)
:=
\frac1n
\sum_{i=-n}^{-1}
\mathbf 1\left\{
X_i^{i+k-2}=a_{-k+1}^{-1},
\,Y_i^{i+k-1}=b_{-k+1}^{0}
\right\}.
\]
By stationarity,
\[
E_\mu\!\left[
\overline\mu_{n;k-1,k}\left(a_{-k+1}^{-1},b_{-k+1}^{0}\right)
\right]
=
\mu_{k-1,k}\left(a_{-k+1}^{-1},b_{-k+1}^{0}\right).
\]
Moreover, $\widehat\mu_{n;k-1,k}$ and $\overline\mu_{n;k-1,k}$ may differ only at the last $k-1$ starting positions. Therefore,
\begin{equation}
\left\|
E_\mu\!\left[\widehat\mu_{n;k-1,k}\right]
-
\mu_{k-1,k}
\right\|_{\rm TV}
\leq
\frac{k-1}{n}.
\label{eq:YX-periodic-boundary-TV}
\end{equation}
Since, by Assumption~{\rm (I)},
$0\leq -\log g_k^{Y|X}\left(b_0\mid a_{-k+1}^{-1},b_{-k+1}^{-1}\right)
\leq \log\frac1{\epsilon_0},$
it follows from \eqref{eq:YX-periodic-boundary-TV} that
\begin{equation}
\left|R_{n,k}^{Y|X}\right|
\leq
\frac{k-1}{n}\log\frac1{\epsilon_0}.
\label{eq:YX-periodic-boundary-error}
\end{equation}
Consequently, combining \eqref{eq:BYX-boundary-KL} and
\eqref{eq:YX-periodic-boundary-error},
\begin{equation}
B_{n,k}^{Y|X}
\leq
\frac{k-1}{n}\log\frac1{\epsilon_0}
+
\left|
E_\mu\!\left[\Delta\widehat D_{n,k}^{Y|X}\right]
\right|.
\label{eq:BYX-KL-reduction}
\end{equation}

Since $\widehat\mu_{n;k-1,k-1}$ and $\mu_{k-1,k-1}$ are, respectively, the marginals of $\widehat\mu_{n;k-1,k}$ and $\mu_{k-1,k}$ obtained by summing over $b_0\in\mathcal Y$, the data-processing inequality for relative entropy yields $D\left(\widehat\mu_{n;k-1,k-1}\Vert\mu_{k-1,k-1}\right) \leq D\left(\widehat\mu_{n;k-1,k}\Vert\mu_{k-1,k}\right)$ \citep{Csiszar1967}. Thus, $\Delta\widehat D_{n,k}^{Y|X}\leq0$, and therefore
\begin{equation}
\left|E_\mu\!\left[\Delta\widehat D_{n,k}^{Y|X}\right]\right|\leq E_\mu\!\left[D\left(
\widehat\mu_{n;k-1,k}\Vert\mu_{k-1,k}\right)\right].
\label{eq:BYX-data-processing}
\end{equation}

We now separate the blocks entirely contained in the observed sample from those completed using the periodic extension. Define
\begin{equation}
\label{eq:auxiliary_v}
\widehat\nu_{n;k-1,k}\left(a_{-k+1}^{-1},b_{-k+1}^{0}\right) := \frac{1}{n-k+1}\sum_{i=-n}^{-k}\mathbf 1\left\{X_i^{i+k-2}=a_{-k+1}^{-1},\,Y_i^{i+k-1}=b_{-k+1}^{0}\right\}.
\end{equation}
Let $\widehat\rho_{n;k-1,k}$ denote the empirical distribution of the remaining $k-1$ blocks completed using the periodic extension. Then
\begin{equation}
\widehat\mu_{n;k-1,k} = \frac{n-k+1}{n}\widehat\nu_{n;k-1,k} + \frac{k-1}{n}\widehat\rho_{n;k-1,k}.
\label{eq:BYX-periodic-decomposition}
\end{equation}
Using \eqref{eq:BYX-periodic-decomposition} and the convexity of relative entropy,
\[
D\left(\widehat\mu_{n;k-1,k}\Vert\mu_{k-1,k}\right) \leq \frac{n-k+1}{n} D\left(\widehat\nu_{n;k-1,k}\Vert\mu_{k-1,k}\right) + \frac{k-1}{n} D\left(\widehat\rho_{n;k-1,k}\Vert\mu_{k-1,k}\right).
\]

By Assumption~{\rm (I)}, every joint $k$-block has probability at least $\epsilon_0^k$. Since $\mu_{k-1,k}$ is a marginal of $\mu_{k,k}$, it follows that $\mu_{k-1,k}\left(a_{-k+1}^{-1},b_{-k+1}^{0}\right)\geq\epsilon_0^k.$ Consequently, for any probability distribution $\rho$ on $\mathcal X^{k-1}\times\mathcal Y^k$,
\[
D\left(\rho\Vert\mu_{k-1,k}\right)\leq k\log\frac1{\epsilon_0}.\]
Hence,
\begin{equation}
E_\mu\!\left[D\left(\widehat\mu_{n;k-1,k}\Vert\mu_{k-1,k}\right)\right] \leq
\frac{n-k+1}{n}E_\mu\!\left[D\left(\widehat\nu_{n;k-1,k}\Vert\mu_{k-1,k}\right)\right] +
\frac{k(k-1)}{n}\log\frac1{\epsilon_0}.
\label{eq:BYX-KL-periodic-bound}
\end{equation}

It remains to control the first term on the right-hand side of \eqref{eq:BYX-KL-periodic-bound}. For $u=\left(a_{-k+1}^{-1},b_{-k+1}^{0}\right)\in\mathcal X^{k-1}\times\mathcal Y^k,$
define
\[
I_i^u := \mathbf 1\left\{X_i^{i+k-2}=a_{-k+1}^{-1},\,Y_i^{i+k-1}=b_{-k+1}^{0}\right\}.
\]
Then, we can rewrite \eqref{eq:auxiliary_v} as
\[
\widehat\nu_{n;k-1,k}(u) = \frac1{n-k+1}\sum_{i=-n}^{-k} I_i^u,
\]
and, by stationarity,
\[
E_\mu\left[\widehat\nu_{n;k-1,k}(u)\right] = \mu_{k-1,k}(u).
\]

Using the standard bound of relative entropy by the $\chi^2$-divergence \citep[Lemma 2.7]{Tsybakov2009},
\[
D\left(\widehat\nu_{n;k-1,k}\Vert\mu_{k-1,k}\right) \leq \sum_{u\in\mathcal X^{k-1}\times\mathcal Y^k}\frac{\left(\widehat\nu_{n;k-1,k}(u)-\mu_{k-1,k}(u)\right)^2}{\mu_{k-1,k}(u)}.
\]
Hence,
\[
E_\mu\left[D\left(\widehat\nu_{n;k-1,k}\Vert\mu_{k-1,k}\right)\right] \leq \sum_{u\in\mathcal X^{k-1}\times\mathcal Y^k}\frac{\operatorname{Var}_\mu\left(\widehat\nu_{n;k-1,k}(u)\right)}{
\mu_{k-1,k}(u)}.
\]

For each $u$, the variable $I_i^u$ is measurable with respect to $\sigma(Z_i,\ldots,Z_{i+k-1})$. Therefore, exactly as in the proof of Theorem~\ref{thm:concentration_entropy_rate}, for $r\geq k$,
\[\left|\operatorname{Cov}_\mu(I_0^u,I_r^u)\right| \leq \min\left\{\mu_{k-1,k}(u),\,\beta_Z(r-k+1)\right\}.\]
For $1\leq r<k$, we simply use
\[
\left|\operatorname{Cov}_\mu(I_0^u,I_r^u)\right|\leq\mu_{k-1,k}(u).
\]
Consequently,
\[
\operatorname{Var}_\mu\left(\widehat\nu_{n;k-1,k}(u)\right) \leq \frac{\mu_{k-1,k}(u)}{n-k+1}\left[2k-1 + 2\sum_{s=1}^{\infty}\min\left\{1,\frac{\beta_Z(s)}{\mu_{k-1,k}(u)}\right\}\right].
\]

By Assumption~{\rm (III)}, $\beta_Z(s)\leq\gamma e^{-\delta s}$, and the same estimate used in the proof of Theorem~\ref{thm:concentration_entropy_rate} gives
\[
\sum_{s=1}^{\infty} \min\left\{1,\frac{\beta_Z(s)}{\mu_{k-1,k}(u)}\right\} \leq \frac{
\log_+\left(\gamma/\mu_{k-1,k}(u)\right)}{\delta} + 1 + \frac{1}{e^\delta-1}.
\]
Moreover, Assumption~{\rm (I)} implies $\mu_{k-1,k}(u)\geq\epsilon_0^k.$ Thus,
\[
\log_+\left(\frac{\gamma}{\mu_{k-1,k}(u)}\right) \leq \log_+\gamma + k\log\frac1{\epsilon_0}.
\]
Recalling $c_1 = 2+\frac{2}{\delta}\log\frac1{\epsilon_0}$ and $c_0 = 1+\frac{2\log_+\gamma}{\delta} +\frac{2}{e^\delta-1},$ we conclude that
\[
\operatorname{Var}_\mu\left(\widehat\nu_{n;k-1,k}(u)\right) \leq \frac{\mu_{k-1,k}(u)}{n-k+1}(c_1k+c_0).
\]
Since there are $|\mathcal X|^{k-1}|\mathcal Y|^k$ possible values of $u$, it follows that
\begin{equation}
E_\mu\left[D\left(\widehat\nu_{n;k-1,k}\Vert\mu_{k-1,k}\right)\right] \leq
\frac{|\mathcal X|^{k-1}|\mathcal Y|^k}{n-k+1}(c_1k+c_0).
\label{eq:BYX-interior-KL-bound}
\end{equation}
Combining \eqref{eq:BYX-interior-KL-bound} with \eqref{eq:BYX-KL-periodic-bound}, we obtain
\begin{equation}
E_\mu\left[
D\left(
\widehat\mu_{n;k-1,k}\Vert\mu_{k-1,k}
\right)
\right]
\leq
\frac{
|\mathcal X|^{k-1}|\mathcal Y|^k
}{n}
(c_1k+c_0)
+
\frac{k(k-1)}{n}\log\frac1{\epsilon_0}.
\label{eq:BYX-full-KL-bound}
\end{equation}
It follows from \eqref{eq:BYX-data-processing} and \eqref{eq:BYX-full-KL-bound} that
\begin{equation}
\left|E_\mu\left[\Delta\widehat D_{n,k}^{Y|X}\right]\right| \leq
\frac{|\mathcal X|^{k-1}|\mathcal Y|^k}{n} (c_1k+c_0) + \frac{k(k-1)}{n}\log\frac1{\epsilon_0}.
\label{eq:BYX-Delta-D-final}
\end{equation}
Finally, substituting \eqref{eq:BYX-Delta-D-final} into \eqref{eq:BYX-KL-reduction} yields
\begin{equation}
B_{n,k}^{Y|X} \leq \frac{|\mathcal X|^{k-1}|\mathcal Y|^k}{n}(c_1k+c_0) + \frac{k^2-1}{n}\log\frac1{\epsilon_0}.
\label{eq:BYX-final-bound}
\end{equation}

Combining \eqref{eq:TE-bias-decomposition}, \eqref{eq:BY-bound}, and \eqref{eq:BYX-final-bound}, we obtain
\begin{align}
\left|E_\mu\!\left[\widehat T_n^k\right] - T(\boldsymbol X\to\boldsymbol Y) \right|
\leq\;& \gamma_k^Y+\gamma_k^{Y|X} \nonumber\\
&+\frac{|\mathcal Y|^k+|\mathcal X|^{k-1}|\mathcal Y|^k}{n}(c_1k+c_0) + \frac{2(k^2-1)}{n}\log\frac1{\epsilon_0}.
\label{eq:TE-bias-combined}
\end{align}
Since $\mathcal Z=\mathcal X\times\mathcal Y$, $|\mathcal Y|^k\leq|\mathcal Z|^k$ and $|\mathcal X|^{k-1}|\mathcal Y|^k \leq |\mathcal Z|^k.$ Thus, \eqref{eq:TE-bias-combined} yields
\begin{equation}
\left|E_\mu\!\left[\widehat T_n^k\right] - T(\boldsymbol X\to\boldsymbol Y) \right| \leq
\gamma_k^Y+\gamma_k^{Y|X} + \frac{2|\mathcal Z|^k}{n}(c_1k+c_0) + \frac{2(k^2-1)}{n}\log\frac1{\epsilon_0}.
\label{eq:TE-bias-final}
\end{equation}

We now take $k=k(n)$ and define
\begin{equation}
r_n := \gamma_{k(n)}^Y + \gamma_{k(n)}^{Y|X} + \frac{2|\mathcal Z|^{k(n)}}{n} \bigl(c_1k(n)+c_0\bigr) + \frac{2(k(n)^2-1)}{n}\log\frac1{\epsilon_0}.
\label{eq:TE-bias-sequence}
\end{equation}
By \eqref{eq:TE-bias-final} and \eqref{eq:TE-bias-sequence},
\begin{equation}
\left|E_\mu\!\left[\widehat T_n^{k(n)}\right] - T(\boldsymbol X\to\boldsymbol Y) \right| \leq r_n.
\label{eq:TE-bias-rn}
\end{equation}
Since $k(n)\stackrel{n \rightarrow \infty}{\longrightarrow}\infty$, the definitions of the entropy rate and the conditional entropy rate imply
\[
\gamma_{k(n)}^Y\stackrel{n \rightarrow \infty}{\longrightarrow} 0
\qquad\text{and}\qquad
\gamma_{k(n)}^{Y|X}\stackrel{n \rightarrow \infty}{\longrightarrow} 0.
\]
Moreover, assume that, for all sufficiently large $n$, $k(n)\leq C\log n,$ $0<C<\frac{1}{\log|\mathcal Z|},$ then by \eqref{eq:TE-bias-sequence}, $r_n\stackrel{n \rightarrow \infty}{\longrightarrow}0.$ By the triangle inequality and \eqref{eq:TE-bias-rn},
\[\left|\widehat T_n^{k(n)} - T(\boldsymbol X\to\boldsymbol Y) \right| \leq \left|\widehat T_n^{k(n)} - E_\mu\!\left[\widehat T_n^{k(n)}\right]\right| + r_n.
\]

Therefore, for every $\xi>0$,
\begin{equation}
\mu\left(\left|\widehat T_n^{k(n)}-T(\boldsymbol X\to\boldsymbol Y)\right|>\xi+r_n\right)
\leq \mu\left(\left|\widehat T_n^{k(n)} - E_\mu\!\left[\widehat T_n^{k(n)}\right]\right|
> \xi\right).
\label{eq:TE-rate-to-mean-reduction}
\end{equation}
By Theorem~\ref{thm:concentration_transfer_entropy_mean}, the right-hand side of
\eqref{eq:TE-rate-to-mean-reduction} is bounded by
\[4\exp\left\{-\frac{2\xi^2 n}{9(1+E[\theta])^2(4k(n)-2)^2(\log n+1)^2}
\right\}.
\]
Since, for all sufficiently large $n$, $k(n)\leq C\log n,$ we have
\begin{equation}
\mu\left(
\left|
\widehat T_n^{k(n)}
-
T(\boldsymbol X\to\boldsymbol Y)
\right|
>
\xi+r_n
\right)
\leq
4\exp\left\{
-\frac{
\xi^2 n
}{
288C^2(1+E[\theta])^2\log^4 n
}
\right\},
\label{eq:TE-rate-concentration-final}
\end{equation}
which completes the proof.

\qed

\subsection{Proof of Corollary \ref{cor:concentration_transfer_entropy_rate}}
\label{sec:proof_corollary_transfer_entropy_rate}

Since $r_n\to0$ and
\[
\sum_{n\ge2}
\exp\left\{
-c\,\frac{n}{\log^4 n}
\right\}
<\infty,
\]
with $c > 0$, Theorem~\ref{thm:concentration_transfer_entropy_rate} and the
Borel--Cantelli lemma \cite[Theorem~2.3.1]{Durrett2019} imply
\[
\widehat T_n^{k(n)}
\left(
X_{-k(n)+1}^{0}\to Y_{-k(n)+1}^{0}
\right)
\longrightarrow
T(\boldsymbol X\to\boldsymbol Y)
\qquad
\mu\text{-a.s.}
\]
\qed

\section{Final remarks}
\label{sec:final_remarks}

We established non-asymptotic concentration inequalities for plug-in estimators of the entropy rate and the transfer entropy rate for stationary finite-alphabet chains with unbounded memory governed by $g$-functions that need not be globally continuous. The coupling-from-the-past assumption controls the fluctuations of the estimators around their expectations, while uniform non-nullness and exponential $\beta$-mixing provide the additional control needed to obtain concentration around the corresponding rates and almost-sure consistency. Natural directions for further work include weakening the mixing and
non-nullness assumptions, as well as establishing central limit theorems and large deviation results under the same non-regular framework.

\section*{Acknowledgements}

\noindent R.F.F. was partially supported by the São Paulo Research Foundation (FAPESP), Brazil, Process Number \#2025/13627-9.

%\appendix
%\section{Example Appendix Section}
%\label{app1}
%%       \citet{<label>} ==> Lamport (1994)
%%       \citep{<label>} ==> (Lamport, 1994)
%%
%Example citation, See \citet{lamport94}.



%\begin{thebibliography}{00}


%\bibitem[Lamport(1994)]{lamport94}
 % Leslie Lamport,
%  \textit{\LaTeX: a document preparation system},
%  Addison Wesley, Massachusetts,
 % 2nd edition,
  %1994.

%\end{thebibliography}
%\end{document}

%\endinput
%%
%% End of file `elsarticle-template-harv.tex'.

\bibliographystyle{elsarticle-harv} 
\bibliography{references}

\end{document}